\documentclass[12pt]{article}
\usepackage{amsfonts, amssymb, amsmath, amsthm, epsfig,cite}

\makeatletter 

\newtheorem{theorem}{Theorem}
\@addtoreset{theorem}{section}

\newtheorem{lemma}{Lemma}
\@addtoreset{lemma}{section}

\@addtoreset{definition}{section}

\newtheorem{remark}{Remark}
\@addtoreset{remark}{section}

\@addtoreset{example}{section}

\newtheorem{condition}{Condition}
\@addtoreset{condition}{section}

\@addtoreset{corollary}{section}

\@addtoreset{figure}{section}

\@addtoreset{equation}{section}

\makeatother

\title{
Asymptotics of solutions for nonlocal elliptic
problems in plane bounded domains
\thanks{This work has been supported by Russian Foundation for Basic
Research (grant~02-01-00312) and by INTAS (grant~YSF~2002-008).} }

\author{Pavel Gurevich}
\date{}

\begin{document}

\maketitle

\begin{abstract}

The paper is devoted to the study of asymptotic behavior of
solutions for nonlocal elliptic problems in weighted spaces. We
deal with the most difficult case when the support of nonlocal
terms intersects with boundary of a plane bounded domain. In this
situation, a general form of the asymptotics is investigated, and
coefficients in the asymptotics are calculated.

\end{abstract}





\bigskip

\begin{center}
 \bf Contents
\end{center}

\ref{sectIntrod}. Introduction \hfill \pageref{sectIntrod}

\ref{sectStatement}. Statement of the problem in a bounded domain
\hfill \pageref{sectStatement}

\ref{sectAsymp}. Asymptotics of solutions for nonlocal problems
\hfill \pageref{sectAsymp}

\ref{sectIndex}. Index of nonlocal problems \hfill
\pageref{sectIndex}

\ref{sectAsymp*}. Asymptotics of solutions for adjoint nonlocal
problems \hfill \pageref{sectAsymp*}

\ref{sectCoef}. Calculation of the coefficients in the asymptotics
formulas \hfill \pageref{sectCoef}

\ref{sectExample}. Example \hfill \pageref{sectExample}

References \hfill \pageref{sectBibliogr}

\bigskip

\section{Introduction}
\phantom{a}~\newline \nopagebreak \label{sectIntrod}

{\bf I.} This work is devoted to the investigation of asymptotic
behavior of solutions for nonlocal elliptic problems. Recently
many mathematicians have been studying nonlocal problems. This
interest is explained, on the one hand, by a significant
theoretical progress in the area and, on the other hand, by a
number of important applications arising in plasma
theory~\cite{BitzSam}, biophysics, theory of diffusion
processes~\cite{Feller,Ventz,SkBook}, modern aircraft technology
(particularly, in the theory of sandwich shells and
plates~\cite{SkBook}), etc.

In the 1-dimensional case the first ones who studied nonlocal
problems were A.~Sommerfeld~\cite{Sommerfeld},
Ya.D.~Tamarkin~\cite{Tamarkin}, M.~Picone~\cite{Picone}. In the
2-dimensional case the earliest paper devoted to nonlocal problems
is due to T.~Carleman~\cite{Carleman}. T.~Carleman searched for a
harmonic function $u$ in a plane bounded domain $G$ subject to a
nonlocal condition connecting the values of the unknown function
in different points of boundary:
$u(x)+bu\big(\omega(x)\big)=g(x)$. Here $\omega:\partial G\to
\partial G$ is a nondegenerate transformation subject to the restriction
$\omega(\omega(x))\equiv x$ (being referred  to as Carleman's
condition in the present time). Such a statement of nonlocal
problems has originated further research into the area of elliptic
problems with nonlocal transformations mapping a boundary onto
itself and with abstract boundary
conditions~\cite{Vishik,Browder,Beals,Antonev}.

In 1969, A.V.~Bitsadze and A.A.~Samarskii~\cite{BitzSam}
considered the following nonlocal problem arising in plasma
theory: find a harmonic function $u(y_1,\ y_2)$ in the rectangle
$G=\{y\in{\mathbb R}^2: -1<y_1<1,\ 0<y_2<1\}$ such that it is
continuous in $\bar G$ and satisfies the conditions
$$
 \begin{array}{c}
  u(y_1,\ 0)=f_1(y_1),\ u(y_1,\ 1)=f_2(y_1),\ -1<y_1<1,\\
  u(-1,\ y_2)=f_3(y_2),\ u(1,\ y_2)=u(0,\ y_2),\ 0<y_2<1,
 \end{array}
$$
where $f_1,\ f_2,\ f_3$ are given continuous functions. We notice
that this problem principally differs from the one studied by
T.~Carleman: now the values of the unknown function on the part of
the boundary $\partial G$ are connected with the values inside the
domain $G$. This problem was solved in~\cite{BitzSam} by reducing
to an integral Fredholm equation and using the maximum principle.
In case of an arbitrary domain and general nonlocal
transformations, it was formulated as an unsolved one.

The most difficult case turns out to deal with the situation when
a part $\Upsilon_1$ of boundary of a domain $G$ is mapped by some
nonlocal transformation $\Omega_1$ on $\Omega_1(\Upsilon_1)$ so
that $\overline{\Omega_1(\Upsilon_1)}\cap\partial
G\ne\varnothing$. Various versions of such problems were
considerer by S.D.~Eidelman and N.V.~Zhitarashu~\cite{ZhEid},
K.Yu.~Kishkis~\cite{Kishk}, A.K.~Gushchin and
V.P.~Mikhailov~\cite{GM}, etc.

\smallskip

Basis of general theory for elliptic equations of order~$2m$ with
general nonlocal conditions was founded by A.L.~Skubachevskii and
his pupils. In a series of works a priori estimates were proved, a
right regularizer was constructed, adjoint problems were studied,
and properties of index in appropriate spaces were established;
spectral properties of some problems were considered~\cite{SkMs83,
SkMs86, SkDu90, SkDu91, SkJMAA,Podiap,SkKov,GurDAN, GurGiess};
asymptotics and smoothness of solutions near some special points
were investigated~\cite{SkMs86, GurPlane}. We remark that the
papers~\cite{SkMs86,SkDu90,SkDu91} were the first ones to deal
with the case $\overline{\Omega_1(\Upsilon_1)}\cap\bar
\Upsilon_1\ne\varnothing$, which had not been previously
considered even for the Laplace equation with nonlocal conditions
in plane domains.

\medskip

{\bf II.} In this paper we investigate the most difficult
situation mentioned above: the support of nonlocal terms can have
a nonempty intersection with boundary of a domain $G$. In that
case, power singularities for solutions near some set ${\cal
K}\subset G$ can appear~\cite{SkMs86, SkRJMP}. Therefore it is
quite natural to study such problems in special weighted spaces
that take into consideration those possible singularities. (The
most convenient spaces turned out to be Kondrat'ev's
ones~\cite{KondrTMMO67}.) Thus we arrive at the question of
asymptotics of solutions near the set ${\cal K}$. In the
paper~\cite{SkMs86}, A.L.~Skubachevskii obtained a general form of
an asymptotics of solutions to problems with nonlocal
transformations coinciding with a rotation operator near the set
${\cal K}$. These theorems were applied to investigation of
smoothness for generalized solutions of nonlocal elliptic problems
(see~\cite{SkMs86,SkRJMP}).

In the present work we generalize the mentioned results of
A.L.~Sku\-ba\-chev\-skii and study the case of arbitrary nonlocal
transformations, linear near the set ${\cal K}$. Simultaneously,
we get a formula connecting the indices of one and the same
nonlocal problem, but being considered in different weighted
spaces.

Moreover, using the results of the paper~\cite{GurPlane} (which
deals with model nonlocal problems in plane angles and in
${\mathbb R}^2\setminus\{0\}$), we get explicit formulas for
calculating the coefficients in the asymptotics of solutions.
These formulas are given both in terms of eigenvectors and
associated vectors of model adjoint problems and in terms of
distributions from the kernel of adjoint problem in a bounded
domain. The latter shows, in particular, that the values of the
coefficients in the asymptotics are the functionals over the
right--hand sides of the nonlocal problem under consideration.
These functionals depend on the data of the problem in the whole
domain, but not only in some neighborhood of the set ${\cal K}$.

We remark that the calculation of the coefficients in the
asymptotics is both important itself and has a direct application
to the question of smoothness of generalized solutions for
nonlocal problems. Roughly speaking, it allows to show that a
generalized solution $u\in W_2^1(G)$ to a nonlocal problem (for an
elliptic 2nd order equation) with a right--hand side $f\in L_2(G)$
is smooth (i.e., $u\in W_2^2(G)$) if and only if the function $f$
satisfies some orthogonality conditions. In a number of cases
these conditions can be verified explicitly.

\smallskip

{\bf III.} The paper is organized as follows. The statement of the
problem and some assumptions concerning nonlocal transformations
are given in section~\ref{sectStatement}. Most of the assumptions
are due to simplify computations throughout the paper. In
section~\ref{sectAsymp} we derive an asymptotics (with yet unknown
coefficients) for solutions to nonlocal problems. Using the
results of section~\ref{sectAsymp}, in section~\ref{sectIndex} we
establish a connection between the indices of one and the same
problem but being considered in different weighted spaces. In
section~\ref{sectAsymp*} we obtain an asymptotics of solutions for
adjoint nonlocal problems. This allows to get in
section~\ref{sectCoef} explicit formulas for the coefficients in
asymptotics of solutions to the original nonlocal problem. In
section~\ref{sectExample} we consider an example illustrating the
results of sections~\ref{sectStatement}--\ref{sectCoef}.

\section{Statement of the problem in a bounded domain}
\phantom{a}~\newline \nopagebreak\label{sectStatement}

Let $G\in{\mathbb R}^2$ be a bounded domain with a boundary
$\partial G=\bigcup\limits_{\sigma=1,2} \bar\Upsilon_\sigma$,
where $\Upsilon_\sigma$ are open (in the topology of $\partial G$)
curves of $C^\infty$ class such that
$\Upsilon_1\cap\Upsilon_2=\varnothing$,
$\bar\Upsilon_1\cap\bar\Upsilon_2=\{g_1,\ h_1\}$. We suppose that
in some neighborhoods of the points $g_1$ and $h_1$ the domain $G$
coincides with an angle.

We denote by ${\bf P}(y,\ D_y)$, $B_{\sigma\mu}(y,\ D_y)$,
$T_{\sigma\mu}(y,\ D_y)$ differential operators of orders $2m$,
$m_{\sigma\mu}$, $m_{\sigma\mu}$ respectively with complex--valued
coefficients from $C^\infty({\mathbb R}^2)$
($m_{\sigma\mu}\le2m-1$, $\sigma=1,\ 2;$ $\mu=1,\ \dots,\ m$). Put
also $B_{\sigma}(y,\ D_y)=\{B_{\sigma\mu}(y,\ D_y)\}_{\mu=1}^m$,
$T_{\sigma}(y,\ D_y)=\{T_{\sigma\mu}(y,\ D_y)\}_{\mu=1}^m$.

\begin{figure}
  \begin{center}
    \leavevmode
    \psfig{file=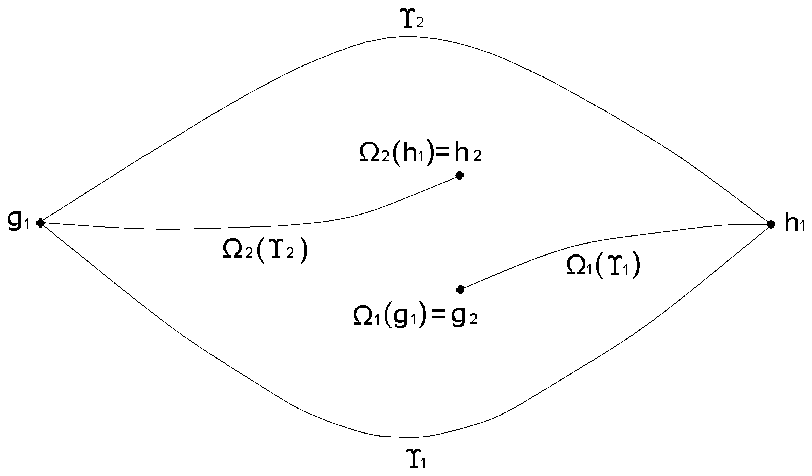}
    \caption{The domain $G$.}
    \label{figDomG}
  \end{center}
\end{figure}

Let $\Omega_{\sigma }$ ($\sigma=1,\ 2$) be an infinitely
differentiable nondegenerate transformation mapping some
neighborhood ${\cal O}_\sigma$ of $\Upsilon_\sigma$ onto
$\Omega_{\sigma }({\cal O}_\sigma)$ such that $\Omega_{\sigma
}(\Upsilon_\sigma)\subset G.$ For definiteness, we consider the
case when $\Omega_1(g_1)=g_2\in G$, $\Omega_2(g_1)=g_1$,
$\Omega_1(h_1)=h_1$, $\Omega_2(h_1)=h_2\in G$. In this work we
also assume that $g_2\notin\overline{\Omega_2(\Upsilon_2)}$,
$h_2\notin\overline{\Omega_1(\Upsilon_1)}$. The last assumption is
made in order to simplify further computations\footnote{If, say,
$g_2\in\overline{\Omega_2(\Upsilon_2)}$, then either $g_2=h_2$ (in
that case, an asymptotics of a solution near the point $g_2$ will
influence not only an asymptotics near $g_1$ but near $h_1$ as
well) or $g_2\in\Omega_2(\Upsilon_2)$ (in that case, one must
study an asymptotics at the additional point
$\Omega_2^{-1}(g_2)\in\Upsilon_2$).}. But, following~\cite{SkMs86,
SkDu91}, we demand (and it is on principle) that following
condition holds:

\begin{condition}
The curves $\Omega_1(\Upsilon_1)$ and $\Omega_2(\Upsilon_2)$ are
not tangent to the boundary $\partial G$ at the ``consistent
points'' $h_1$ and $g_1$ respectively (see Fig.~\ref{figDomG}).
\end{condition}

We also suppose for simplicity that the transformations
$\Omega_{\sigma }(y)$ are linear near the points $g_1$ and $h_1$.

We introduce the set ${\cal K}=\{g_1,\ h_1,\ g_2,\ h_2\}$ and
consider the nonlocal elliptic problem
\begin{gather}
 {\bf P}(y,\ D_y) u=f(y) \quad (y\in G\backslash{\cal K}),\label{eqPinG}\\
 \begin{aligned}\label{eqBinG}
{\bf B}_{\sigma}(y,\ D_y)u\equiv B_{\sigma}(y,\
D_y)u|_{\Upsilon_\sigma}+
\big(T_{\sigma}(y,\ D_y)u\big)\big(\Omega_{\sigma}(y)\big)|_{\Upsilon_\sigma}=
f_{\sigma}(y)\\
    (y\in \Upsilon_\sigma;\ \sigma=1,\ 2).
  \end{aligned}
\end{gather}
Here $\big(T_{\sigma}(y,\
D_y)u\big)\big(\Omega_{\sigma}(y)\big)=T_{\sigma}(y',\
D_{y'})u(y')|_{y'=\Omega_{\sigma}(y)}$;
$f_{\sigma}=\{f_{\sigma\mu}\}_{\mu=1}^m$.

\begin{remark}
 The results of this paper are generalized for the case when the
boundary $\partial G$ consists of a finite number of smooth curves
$\Upsilon_\sigma$, $\sigma=1,\ \dots,\ N$, and nonlocal conditions
on each $\Upsilon_\sigma$ contain a finite number of nonlocal
terms with different transformations. Moreover, these
transformations can map ``the consistent points'' (which are $g_1$
and $h_1$ in our case) both to the boundary $\partial G$ and
inside the domain $G$, forming finite orbits.
\end{remark}

We introduce the space $H_a^l(G)$ as a completion of the set
$C_0^\infty(\bar G\backslash {\cal K})$ in the norm
$$
\|u\|_{H_a^l(G)}=\left(\sum\limits_{|\alpha|\le l}\int\limits_G
\rho^{2(a-l+|\alpha|)}|D^\alpha u|^2 dy\right)^{1/2}.
$$
Here $C_0^\infty(\bar G\backslash {\cal K})$ is the set of
infinitely differentiable functions with compact supports
contained in $\bar G\backslash {\cal K}$; $l\ge 0$ is an integer;
$a\in {\mathbb R};$ $\rho=\rho(y)={\rm dist}(y,\ {\cal K})$.

If, instead of the domain $G$, one considers an angle with a
vertex $g$ or some neighborhood of a point $g$, then one must put
${\cal K}=\{g\}$ in the definition of the weighted space.

By $H_a^{l-1/2}(\Upsilon)$ we denote the space of traces on a
smooth curve $\Upsilon\subset\bar G$ with the norm $
\|\psi\|_{H_a^{l-1/2}(\Upsilon)}=\inf\|u\|_{H_a^l(G)}\quad (u\in
H_a^l(G):\ u|_\Upsilon=\psi). $

Let us introduce the operator
\begin{multline}\notag
 {\bf L}=\{{\bf P}(y,\ D_y),\ {\bf B}_{\sigma}(y,\ D_y)\}: \\
 H_a^{l+2m}(G)\to  H_a^l(G,\ \Upsilon)\stackrel{def}{=}
H_a^l(G)\times\prod\limits_{\sigma=1,2}\prod\limits_{\mu=1}^m
H_a^{l+2m-m_{\sigma\mu}-1/2}(\Upsilon_\sigma),
\end{multline}
which corresponds to nonlocal problem~(\ref{eqPinG}),
(\ref{eqBinG}).

\medskip

Throughout the paper we assume that the operators ${\bf P}(y,\
D_y)$ and $B_{\sigma}(y,\ D_y)$ satisfy the following conditions
(see, e.g.,~\cite[Chapter~2, section~1]{LM}).

\begin{condition}\label{condEllipPinG}
For all $y\in\bar G$ the operator ${\bf P}(y,\ D_y)$ is properly
elliptic.
\end{condition}
\begin{condition}\label{condComplBinG}
For $\sigma=1,\ 2$ and $y\in\bar\Upsilon_\sigma$ the system
$B_{\sigma}(y,\ D_y)=\{B_{\sigma\mu}(y,\ D_y)\}_{\mu=1}^m$ covers
the operator ${\bf P}(y,\ D_y)$.
\end{condition}

Remark that we do not impose any restrictions on the nonlocal
operators $T_{\sigma\mu}(y,\ D_y)$ but the natural restriction on
their orders.

\section{Asymptotics of solutions for nonlocal problems}
\phantom{a}~\newline \nopagebreak\label{sectAsymp}

{\bf I.} In this section we obtain an asymptotics of a given
solution $u\in H_a^{l+2m}(G)$ for problem~(\ref{eqPinG}),
(\ref{eqBinG}) with a right--hand side $\{f,\ f_{\sigma}\}\in
H_{a_1}^l(G,\ \Upsilon)$, $0<a-a_1<1$.

Notice that the violation of the inequality $a-a_1<1$ means that
$\{f,\ f_{\sigma}\}$ is ``too regular''. In that case, exact
results should yield more terms in asymptotics in comparison with
our case. This situation can be investigated in the way similar
to~\cite{SkMs86} (see also~\cite{KondrTMMO67, NP}). Namely, one
should consider corresponding equations for a residue in the
asymptotics formula and apply to them the results obtained for the
case $a-a_1<1$. We are not going to do this here since detailed
computations would lead to enormous enlargement of the paper,
giving no essentially new results (with respect to the present
work and to the paper~\cite{SkMs86}). The same remark is valid for
the case when $u\in H_a^{l+2m}(G)$, $\{f,\ f_{\sigma}\}\in
H_{a_1}^{l_1}(G,\ \Upsilon)$, and $l_1\ne l$.

\smallskip

The asymptotics will be found with the help of eigenvalues and
corresponding Jordan chains of some holomorphic operator--valued
functions. Therefore let us remind some relevant definitions and
facts (see~\cite{GS}).

Suppose $\tilde{\cal L}(\lambda):H_1\to H_2$ is a holomorphic
operator--valued function, $H_1,\ H_2$ are Hilbert spaces. A
holomorphic at a point $\lambda_0$
vector--function~$\varphi(\lambda)$ with the values in~$H_1$ is
called {\it a root function} of the operator~$\tilde{\cal
L}(\lambda)$ at $\lambda_0$ if $\varphi(\lambda_0)\ne0$ and the
vector--function~$\tilde{\cal L}(\lambda)\varphi(\lambda)$ is
equal to 0 at~$\lambda_0$. If~$\tilde{\cal L}(\lambda)$ has at
least one root function at the point~$\lambda_0$, then $\lambda_0$
is called {\it an eigenvalue} of $\tilde{\cal L}(\lambda)$.
Multiplicity of zero for the vector--function~$\tilde{\cal
L}(\lambda)\varphi(\lambda)$ at the point~$\lambda_0$ is called
{\it a multiplicity of the root function} $\varphi(\lambda)$; the
vector~$\varphi^{(0)}=\varphi(\lambda_0)$ is called {\it an
eigenvector} corresponding to the eigenvalue~$\lambda_0$. Let
$\varphi(\lambda)$ be a root function at the point~$\lambda_0$ of
multiplicity~$\varkappa$, and
$\varphi(\lambda)=\sum\limits_{j=0}^{\infty}(\lambda-\lambda_0)^j\varphi^{(j)}$.
Then the vectors $\varphi^{(1)},\ \dots,\ \varphi^{(\varkappa-1)}$
are called {\it associated with the eigenvector} $\varphi_0$, and
the ordered set~$\varphi^{(0)},\ \dots,\ \varphi^{(\varkappa-1)}$
is called  {\it a Jordan chain} corresponding to the
eigenvalue~$\lambda_0$. {\it Rank} of the
eigenvector~$\varphi^{(0)}$ (${\rm rank\,}\varphi^{(0)}$) is the
maximum of multiplicities of all root functions such that
$\varphi(\lambda_0)=\varphi^{(0)}$.

Let an eigenvalue $\lambda_0$ of the operator~$\tilde{\cal
L}(\lambda)$ be isolated, ${\rm dim\,}{\rm ker\,}\tilde{\cal
L}(\lambda_0)<\infty$,  and rank of $\lambda_0$ finite. Suppose
$J={\rm dim\,}{\rm ker\,}\tilde{\cal L}(\lambda_0)$ and
$\varphi^{(0,1)},\ \dots,\ \varphi^{(0,J)}$ is a system of
linearly independent eigenvectors such that ${\rm
rank\,}\varphi^{(0,1)}$ is the greatest of ranks of all
eigenvectors corresponding to the eigenvalue~$\lambda_0$, and
${\rm rank\,}\varphi^{(0,j)}$ ($j=2,\ \dots,\ J$) is the greatest
of ranks of eigenvectors from some orthogonal supplement in~${\rm
ker\,}\tilde{\cal L}(\lambda_0)$ to the linear manifold of the
vectors~$\varphi^{(0,1)},\ \dots,\ \varphi^{(0,j-1)}$. The
numbers~$\varkappa_j={\rm rank\,}\varphi^{(0,j)}$ are called {\it
partial multiplicities} of the eigenvalue~$\lambda_0$, and the
sum~$\varkappa_1+\dots+\varkappa_J$ is called {\it a (full)
multiplicity} of $\lambda_0$. If the vectors~$\varphi^{(0,j)},\
\dots,\ \varphi^{(\varkappa_j-1,j)}$ form a Jordan chain for
every~$j=1,\ \dots,\ J$, then the set of vectors
$\{\varphi^{(0,j)},\ \dots,\ \varphi^{(\varkappa_j-1,j)}:j=1,\
\dots,\ J\}$ is called {\it a canonical system of Jordan chains}
corresponding to the eigenvalue~$\lambda_0$.

\medskip

{\bf II.} At first let us consider an asymptotics of the solution
$u$ for problem~(\ref{eqPinG}), (\ref{eqBinG}) near the point
$g_2$. In this case we will see that the asymptotics is defined by
a model ``local'' problem in ${\mathbb R}^2\setminus\{g_2\}$ with
a ``regular'' right--hand side. Such a problem was studied
in~\cite[section~5]{GurPlane}. Thereafter we will consider the
asymptotics near the point $g_1$. In that case, we will arrive at
a model nonlocal problem in some angle~$K$ with a right--hand side
being a sum of ``regular'' and ``special'' functions. The
asymptotics of the ``special'' one will be defined by the
asymptotics of the solution $u$ near $g_2$, which is explained by
the presence of the nonlocal transformation $\Omega_1$. Then the
results of~\cite{GurPlane} will be applied to this model problem.

Thus we fix a neighborhood ${\cal V}(g_2)$ of $g_2$ such that
$\overline{{\cal V}(g_2)}\cap\partial G=\varnothing$ and
$\overline{{\cal V}(g_2)}\cap\{h_2\}=\varnothing$. One can see
that an asymptotic behavior of $u$ in ${\cal V}(g_2)$ does not
depend on nonlocal conditions~(\ref{eqBinG}), but is defined only
by the equation
\begin{equation}\label{eqP2}
 {\bf P}(y,\ D_y) u=f(y) \quad \big(y\in{\cal V}(g_2)\big).
\end{equation}
Let ${\cal P}(D_y)$ be the principal homogeneous part of the
operator ${\bf P}(g_2,\ D_y)$. Then equation~(\ref{eqP2}) can be
written in the form
\begin{equation}\label{eqPinVg2}
 {\cal P}(D_y)u(y)=\hat f(y)\quad \big(y\in{\cal V}(g_2)\big),
\end{equation}
where $\hat f$, by virtue of the condition $0<a-a_1<1$, belongs to
the space $H_{a_1}^l\big({\cal V}(g_2)\big)$\footnote{To show that
$\hat f\in H_{a_1}^l\big({\cal V}(g_2)\big)$, one must estimate
the expressions of the two types: 1)~$p_\alpha(y)D^\alpha u$,
$|\alpha|\le 2m-1$, and
2)~$\big(p_\alpha(y)-p_\alpha(0)\big)D^\alpha u$, $|\alpha|=2m$,
where $p_\alpha$ are infinitely differentiable coefficients of
${\bf P}(y,\ D_y)$. The 1st one is estimated by direct use of the
condition $0<a-a_1<1$, while the 2nd one needs additional
application of Lemma~3.3$'$~\cite{KondrTMMO67}. Further, in
analogous situations, we will omit these explanations.}. We
introduce the bounded operator
$$
 {\cal L}_2={\cal P}(D_y):H_a^{l+2m}({\mathbb R}^2)\to  H_a^l({\mathbb R}^2),
$$
where, defining the weighted spaces, one must put ${\cal
K}=\{g_2\}$.

We write the operator ${\cal P}(D_y)$ in polar coordinates with
the pole at the point $g_2$: ${\cal P}(D_y)=r^{-2m}\tilde{\cal
P}(\omega,\ D_\omega,\ rD_r)$, where
$D_\omega=-i\frac{\displaystyle\partial}{\displaystyle\partial\omega},\
D_r=-i\frac{\displaystyle\partial}{\displaystyle\partial r}.$

Let us introduce the operator--valued function
$$
 \tilde{\cal L}_2(\lambda)=\tilde{\cal P}(\omega,\ D_\omega,\ \lambda):
 W_{2,2\pi}^{l+2m}(0,\ 2\pi)\to W_{2,2\pi}^l(0,\ 2\pi).
$$
Here $W_{2,2\pi}^{l}(0,\ 2\pi)$ is the closure of the set of
infinitely differentiable $2\pi$-periodic functions in
$W_2^{l}(0,\ 2\pi).$

The operator $\tilde{\cal L}_2(\lambda)$ is obtained from the
operator ${\cal L}_2$ by passing to polar coordinates, followed by
the Mellin transformation with respect to $r$:
$$
 \tilde u(\lambda)=(2\pi)^{-1/2}\int\limits_0^\infty
 r^{-i\lambda-1}u(r)\,dr.
$$

From~\cite[section~1]{SkMs86} it follows that there exists a
finite--meromorphic operator--valued function $\tilde{\cal
R}_2(\lambda)$ such that its poles (except, maybe, a finite number
of them) are located inside a double angle of opening less than
$\pi$ containing the imaginary axis; moreover, if $\lambda$ is not
a pole of $\tilde{\cal R}_2(\lambda)$, then $\tilde{\cal
R}_2(\lambda)$ is inverse to the operator $\tilde{\cal
L}_2(\lambda).$ Thus a number $\lambda$ is a pole of $\tilde{\cal
R}_2(\lambda)$ if and only if $\lambda$ is an eigenvalue (of
finite multiplicity) of $\tilde{\cal L}_2(\lambda)$.

If the line ${\rm Im\,}\lambda=a+1-l-2m$ contains no eigenvalues
of $\tilde{\cal L}_2(\lambda)$, then, by virtue
of~\cite[section~1]{SkMs86}, the operator ${\cal L}_2$ is an
isomorphism.

In order to formulate a theorem concerning an asymptotics near
$g_2$, let us introduce some denotation. Suppose $\lambda_2$ is an
eigenvalue of $\tilde{\cal L}_2(\lambda)$,
\begin{equation}\label{eqJordan2}
\{\varphi_{2}^{(0,\zeta)},\ \dots,\
\varphi_{2}^{(\varkappa_{\zeta,2}-1,\zeta)}: \zeta=1,\ \dots,\
J_{2}\}
\end{equation}
is a canonical system of Jordan chains of the operator
$\tilde{\cal L}_2(\lambda)$ corresponding to the eigenvalue
$\lambda_{2}$.

Consider the vector $u_2=\{u_2^{(k,\zeta)}\}$, where
\begin{equation}\label{eqPowerSolLambda_n2}
   u_{2}^{(k,\zeta)}(\omega,\ r)=r^{i\lambda_{2}}\sum\limits_{q=0}^k
   \frac{\displaystyle 1}{\displaystyle q!}(i\ln r)^q\varphi_{2}^{(k-q,\zeta)}(\omega),
\end{equation}
$(\omega,\ r)$ are polar coordinates with the pole at the point
$g_2$.

Notice that (see~\cite[section~5]{GurPlane}) the vector $u_{2}$,
the components $u_{2}^{(k,\zeta)}$ of which are defined
by~(\ref{eqPowerSolLambda_n2}), satisfies the relation
\begin{equation}\label{eqL2u2=0}
 {\cal L}_2 u_{2}=0.
\end{equation}

\begin{theorem}\label{thAsymp2}
Let the lines ${\rm Im\,}\lambda=a_1+1-l-2m$, ${\rm
Im\,}\lambda=a+1-l-2m$ contain no eigenvalues of $\tilde{\cal
L}_2(\lambda)$ and the strip $a_1+1-l-2m<{\rm
Im\,}\lambda<a+1-l-2m$ contain the only eigenvalue $\lambda_{2}$
of $\tilde{\cal L}_2(\lambda)$. Then
 \begin{equation}\label{eqAsymp2}
  u(y)=c_2u_2(y)+\hat u(y)\quad \big(y\in{\cal V}(g_2)\big)\footnotemark{}.
 \end{equation}
 \footnotetext{In formula~(\ref{eqAsymp2}) and further the
expressions such as $c_{2}u_{2}$ are calculated in the following
way: $
 c_{2}u_{2}=
\sum\limits_{\zeta=1}^{J_{2}}\sum\limits_{k=0}^{\varkappa_{\zeta,2}-1}
  c_{2}^{(k,\zeta)}u_{2}^{(k,\zeta)}.
$ } Here $u_2=\{u_2^{(k,\zeta)}\}$, $u_2^{(k,\zeta)}$ are defined
by~(\ref{eqPowerSolLambda_n2}); $c_{2}=\{c_{2}^{(k,\zeta)}\}$ is a
vector of some constants; $\hat u\in H_{a_1}^{l+2m}\big({\cal
V}(g_2)\big)$\footnote{The results of this work are evidently
generalized for the case when the strip $a_1+1-l-2m<{\rm
Im\,}\lambda<a+1-l-2m$ contains a finite number of eigenvalues of
$\tilde{\cal
L}_2(\lambda)$.}\label{fManyLambda}\newcounter{qountManyLambda}
\setcounter{qountManyLambda}{\value{footnote}}.
\end{theorem}
\begin{proof}
Introduce the cut--off function $\eta\in C^\infty({\mathbb R}^2)$
equal to 1 in some neighborhood of the point $g_2$ and vanishing
outside ${\cal V}(g_2)$. Suppose that the function $\eta u$ is
defined in the whole of ${\mathbb R}^2$, being equal to 0 outside
${\cal V}(g_2)$. Then from~(\ref{eqPinVg2}) and Leibnitz's
formula, it follows that
$$
{\cal L}_2(\eta u)\in H_{a_1}^{l+2m}({\mathbb R}^2).
$$
Now it remains only to apply Theorem~5.1~\cite{GurPlane}, which
establishes the asymptotics of solutions for nonlocal problems in
${\mathbb R}^2\setminus\{g_2\}$.
\end{proof}

\begin{remark}\label{rSuperfluous2}
In fact, the assumption that the line ${\rm Im\,}\lambda=a+1-l-2m$
contains no eigenvalues of $\tilde{\cal L}_2(\lambda)$ is
superfluous. Theorem~\ref{thAsymp2} remains valid even if it is
violated (see Remark~5.1~\cite{GurPlane}). But this assumption
will be used for studying the adjoint nonlocal problem and for
calculating the coefficients~$c_{2}^{(k,\zeta)}$. However, this
assumption does not lead to the loss in generality. Indeed, one
can find an $\varepsilon$, $0<\varepsilon<a-a_1$, such that the
strip $a-\varepsilon+1-l-2m\le{\rm Im\,}\lambda\le a+1-l-2m$
contains no eigenvalues of $\tilde{\cal L}_2(\lambda)$, and
therefore (see~\cite[section~5]{GurPlane}) $u\in
H_{a-\varepsilon}^{l+2m}({\cal V}(g_2))$. Hence we arrive at the
situation of Theorem~\ref{thAsymp2}.
\end{remark}

\begin{remark}\label{rNoEigenval2}
From the results of~\cite{GurPlane}, proof of
Theorem~\ref{thAsymp2}, and Remark~\ref{rSuperfluous2}, it follows
that if the strip $a_1+1-l-2m\le{\rm Im\,}\lambda<a+1-l-2m$ has no
eigenvalues of $\tilde{\cal L}_2(\lambda)$, then $u\in
H_{a_1}^{l+2m}\big({\cal V}(g_2)\big)$ for any right--hand side
$\{f,\ f_{\sigma}\}\in H_{a_1}^l(G,\ \Upsilon)$.
\end{remark}

\medskip

{\bf III.} Now we consider an asymptotics of the solution $u$ for
problem~(\ref{eqPinG}), (\ref{eqBinG}) near the point $g_1$. Fix a
neighborhood ${\cal V}(g_1)$ of $g_1$ such that
\begin{equation}\label{eqVg1_1}
\overline{{\cal
V}(g_1)}\cap\overline{\Omega_1(\Upsilon_1)}=\varnothing\quad
\mbox{and}\quad \overline{{\cal V}(g_1)}\cap\{h_2\}=\varnothing.
\end{equation}

Then one can see that an asymptotic behavior of the solution $u$
is defined by the problem
\begin{gather}
 {\bf P}(y,\ D_y) u=f(y) \quad (y\in{\cal V}(g_1)\cap
 G),\label{eqP1}\\
\begin{aligned}\label{eqB1}
B_{1}(y,\ D_y)u|_{{\cal V}(g_1)\cap\Upsilon_1}=f_1(y)-
\big(T_{1}(y,\ D_y)u\big)\big(\Omega_{1}(y)\big)|_{{\cal V}(g_1)\cap\Upsilon_1} \\
 (y\in {\cal V}(g_1)\cap\Upsilon_1),\\
B_{2}(y,\ D_y)u|_{{\cal V}(g_1)\cap\Upsilon_2}+\big(T_{2}(y,\
D_y)u\big)\big(\Omega_{2}(y)\big)|_{{\cal
V}(g_1)\cap\Upsilon_2}=f_2(y)
 \\
  (y\in {\cal V}(g_1)\cap\Upsilon_2).
 \end{aligned}
\end{gather}

Let ${\cal P}(D_y),$ $B_{\sigma}(D_y)$, $T_{2}(D_y)$ be the
principal homogeneous parts of the operators ${\bf P}(g_1,\ D_y),$
$B_{\sigma}(g_1,\ D_y)$, $T_{2}(g_1,\ D_y)$
respectively\footnote{Notice that earlier, in this section, we
denoted by ${\cal P}(D_y)$ the principal homogeneous part of the
operator ${\bf P}(g_2,\ D_y)$. To be strict we had to denote these
operators by different symbols. But we do not do it since
throughout the paper it will always be clear from the context
whether we consider the principal homogeneous part of ${\bf P}(y,\
D_y)$ at $g_1$ or at $g_2$.}. Let $T_{1}(D_y)$ be the principal
homogeneous part of the operator $T_{1}(g_2,\ D_y)$.

From now on we shall suppose that the origin coincides with the
point $g_1$: $g_1=0$, and
\begin{equation}\label{eqVg1_2}
{\cal V}(g_1)\cap G={\cal V}(0)\cap K,
\end{equation}
where $K$ is the plane angle: $
 K=\{y\in{\mathbb R}^2:\ r>0,\ b_1<\omega<b_2\}
$
with the arms $
 \gamma_\sigma=\{y\in{\mathbb R}^2:\ r>0,\ \omega=b_\sigma\},\
 \sigma=1,\ 2.
$ Here $(\omega,\ r)$ are polar coordinates with the pole at the
point $g_1=0$, $0<b_1<b_2<2\pi$.

According to the assumptions of section~\ref{sectStatement}, the
transformations $\Omega_\sigma(y)$ are linear in ${\cal
V}(g_1)={\cal V}(0)$. Let $\Omega_1(y)$ ($y\in{\cal V}(g_1)$) be a
composition of a rotation and an expansion with respect to $g_1$,
and the shift by the vector $\overrightarrow{g_1g_2}$. Let
$\Omega_2(y)$ ($y\in{\cal V}(g_1)$) coincide with the linear
operator ${\cal G}_2$ of a rotation by an angle $\omega_2$
($b_1<b_2+\omega_2<b_2$) and an expansion with a coefficient
$\beta_2>0$.

Let the neighborhood ${\cal V}(g_1)={\cal V}(0)$ be so small that
$ \{\Omega_1(y):\ y\in{\cal V}(g_1)\}\subset {\cal V}(g_2)$ and
relations~(\ref{eqVg1_1}), (\ref{eqVg1_2}) are fulfilled with the
set $\{\Omega_2(y): y\in{\cal V}(g_1)\}$ substituted for ${\cal
V}(g_1)$ (which is related to~(\ref{eqB1})). (Mention that this
requirement is automatically fulfilled whenever the expansion
coefficients for the transformations $\Omega_\sigma(y)$ near the
point $g_1$ are less or equal to 1.)

Now, using asymptotics formula~(\ref{eqAsymp2}) for the solution
$u$ near $g_2$ and Leibnitz's formula, we get that
problem~(\ref{eqP1}), (\ref{eqB1}) in ${\cal V}(g_1)\cap G$ is
equivalent to the following one in ${\cal V}(0)\cap K$:
\begin{gather}
 {\cal P}(D_y) u=\hat f(y) \quad (y\in{\cal V}(0)\cap
 K),\label{eqP10}\\
\begin{aligned}\label{eqB10}
B_{1}(D_y)u|_{{\cal V}(0)\cap\gamma_1}=\hat f_1(y)-c_2f_{12}(y)
\quad (y\in {\cal
V}(0)\cap\gamma_1),\\
B_{2}(D_y)u|_{{\cal V}(0) \cap\gamma_2}+(T_{2}(D_y)u)({\cal
G}_{2}y)|_{{\cal V}(0)\cap\gamma_2}=\hat f_2(y)
 \quad (y\in {\cal V}(0)\cap\gamma_2).
 \end{aligned}
\end{gather}
Here $\hat f=f-\big({\bf P}(y,\ D_y)-{\cal P}(D_y)\big)u$,
\begin{equation}\label{eqf12}
f_{12}=\big(T_1(D_y)u_2\big)\big(\Omega_1(y)\big)|_{{\cal
V}(0)\cap\gamma_1},
\end{equation}
\begin{multline}\notag
\hat f_1=f_1-\big(B_1(y,\ D_y)-B_1(D_y)\big)u|_{{\cal
V}(0)\cap\gamma_1}-\\
\big(T_1(y,\ D_y)-T_1(D_y))u\big)\big(\Omega_1(y)\big)|_{{\cal
V}(0)\cap\gamma_1}-\big(T_1(D_y)\hat
u\big)\big(\Omega_1(y)\big)|_{{\cal V}(0)\cap\gamma_1},
\end{multline}
\begin{multline}\notag
\hat f_2=f_2-\big(B_2(y,\ D_y)-B_2(D_y)\big)u|_{{\cal
V}(0)\cap\gamma_2}-\\
\big((T_2(y,\ D_y)-T_2(D_y))u\big)\big({\cal
G}_2y\big)|_{{\cal V}(0)\cap\gamma_2}.
\end{multline}
Since $T_{1\mu}(D_y)$ is a homogeneous operator of order
$m_{1\mu}$, from~(\ref{eqf12}) and~(\ref{eqPowerSolLambda_n2}) it
follows that the components $f_{12\mu}^{(k,\zeta)}$ of the vector
$f_{12}=\{f_{12\mu}^{(k,\zeta)}\}$ are linear combinations of the
functions $r^{i\lambda_{2}-m_{1\mu}}(i\ln r)^q$, $0\le q\le k$.

Moreover, by virtue of the condition $0<a-a_1<1$, we have $\hat
f\in H_{a_1}^{l}({\cal V}(0)\cap K)$, $\hat f_{\sigma\mu}\in
H_{a_1}^{l+2m-m_{\sigma\mu}-1/2}({\cal V}(0)\cap\gamma_\sigma)$.

Thus we see that~(\ref{eqP10}), (\ref{eqB10}) is a model nonlocal
problem in ${\cal V}(0)\cap K$ with the right--hand side being the
sum of the ``regular'' function $\{\hat f,\ \hat f_1,\ \hat f_2\}$
and the ``special'' function $\{0,\ -c_2f_{12},\ 0\}$. The
asymptotics of the vector $f_{12}$ is defined by the asymptotics
of the solution $u$ near the point $g_2$, i.e., by the vector
$u_2$ (see~(\ref{eqPowerSolLambda_n2})).

Now we are to apply the results of~\cite{GurPlane}. Put
\begin{equation}\label{eqcalB1}
 \begin{array}{cl}
 &{\cal B}_1(D_y)u=B_1(D_y)u|_{\gamma_1},\\
 &{\cal B}_2(D_y)u=B_2(D_y)u|_{\gamma_2}+(T_2(D_y)u)({\cal G}_2y)|_{\gamma_2}
 \end{array}
\end{equation}
and introduce the bounded operator
\begin{multline}\notag
{\cal L}_1=\{{\cal P}(D_y),\ {\cal B}_{\sigma}(D_y)\}:\\
H_a^{l+2m}(K)\to  H_a^{l}(K,\
\gamma)\stackrel{def}{=}H_a^l(K)\times\prod\limits_{\sigma=1,2}\prod\limits_{\mu=1}^m
H_a^{l+2m-m_{\sigma\mu}-1/2}(\gamma_\sigma),
\end{multline}
which corresponds to the model nonlocal problem in the angle $K$.

Write the operators involved into ${\cal L}_1$ in polar
coordinates: ${\cal P}(D_y)=r^{-2m}\tilde{\cal P}(\omega,\
D_\omega,\ rD_r),$ $B_{\sigma}(D_y)=\{r^{-m_{\sigma\mu}}\tilde
B_{\sigma\mu}(\omega,\ D_\omega,\ rD_r)\}_{\mu=1}^m$,
$T_{\sigma}(D_y)=\{r^{-m_{\sigma\mu}}\tilde T_{\sigma\mu}(\omega,\
D_\omega,\ rD_r)\}_{\mu=1}^m$.

Consider the operator--valued function
$$
\tilde{\cal L}_1(\lambda):W_2^{l+2m}(b_1,\ b_2)\to W_2^{l}[b_1,\
b_2]\stackrel{def}{=}W_2^l(b_1,\ b_2) \times{\mathbb C}^{2m}
$$
given by
 \begin{multline}\notag
\tilde{\cal L}_1(\lambda)\varphi=\{\tilde{\cal P}(\omega,\
D_\omega,\ \lambda)\varphi,\ {\tilde B}_{1\mu}(\omega,\ D_\omega,\
\lambda)\varphi(\omega)|_{\omega=b_1},\\
{\tilde B}_{2\mu}(\omega,\ D_\omega,\
\lambda)\varphi(\omega)|_{\omega=b_2}+\beta_2^{i\lambda-m_{2\mu}}
{\tilde T}_{2\mu}(\omega+\omega_2,\ D_\omega,\
\lambda)\varphi(\omega+\omega_2)|_{\omega=b_2}\}.
\end{multline}
The operator $\tilde{\cal L}_1(\lambda)$ is obtained from the
operator ${\cal L}_1$ by passing to polar coordinates, followed by
the Mellin transformation with respect to $r$.

From Lemmas~2.1, 2.2~\cite{SkDu90} it follows that there exists a
finite--mero\-mor\-phic operator--valued function $\tilde{\cal
R}_1(\lambda)$ such that its poles (except, maybe, a finite number
of them) are located inside a double angle of opening less than
$\pi$ containing the imaginary axis; moreover, if $\lambda$ is not
a pole of $\tilde{\cal R}_1(\lambda)$, then $\tilde{\cal
R}_1(\lambda)$ is inverse to the operator $\tilde{\cal
L}_1(\lambda).$ Thus a number $\lambda$ is a pole of $\tilde{\cal
R}_1(\lambda)$ if and only if $\lambda$ is an eigenvalue (of
finite multiplicity) of $\tilde{\cal L}_1(\lambda)$.

If the line ${\rm Im\,}\lambda=a+1-l-2m$ contains no eigenvalues
of $\tilde{\cal L}_1(\lambda)$, then, by virtue of
Theorem~2.1~\cite{SkDu90}, the operator ${\cal L}_1$ is an
isomorphism.

 In order to formulate a
theorem concerning an asymptotics near $g_1$, let us introduce
some denotation. Suppose $\lambda_1$ is an eigenvalue of the
operator $\tilde{\cal L}_1(\lambda)$ located inside the strip
$a_1+1-l-2m<{\rm Im\,}\lambda<a+1-l-2m$,
\begin{equation}\label{eqJordan1}
\{\varphi_{1}^{(0,\zeta)},\ \dots,\
\varphi_{1}^{(\varkappa_{\zeta,1}-1,\zeta)}: \zeta=1,\ \dots,\
J_{1}\}
\end{equation}
is a canonical system of Jordan chains of the operator
$\tilde{\cal L}_1(\lambda)$ corresponding to the eigenvalue
$\lambda_{1}$.

Consider the vector $u_1=\{u_1^{(k,\zeta)}\}$, where
\begin{equation}\label{eqPowerSolLambda_n1}
u_{1}^{(k,\zeta)}(\omega,\ r)=r^{i\lambda_1}\sum\limits_{q=0}^k
\frac{\displaystyle 1}{\displaystyle q!}(i\ln
r)^q\varphi_{1}^{(k-q,\zeta)}(\omega),
\end{equation}
$(\omega,\ r)$ are polar coordinates with the pole at the point
$g_1=0$.

Notice that (see Lemma~2.2~\cite{GurPlane}) the vector $u_{1}$,
the elements $u_{1}^{(k,\zeta)}$ of which are defined
by~(\ref{eqPowerSolLambda_n1}), satisfies the relation
\begin{equation}\label{eqL1u1=0}
 {\cal L}_1 u_{1}=0.
\end{equation}

If $\lambda_2$ is an eigenvalue of $\tilde{\cal L}_1(\lambda)$
(i.e., $\lambda_2=\lambda_1$), then denote by
$\varkappa(\lambda_2)$ the greatest of partial multiplicities of
$\lambda_2$. If $\lambda_2$ is not an eigenvalue of $\tilde{\cal
L}_1(\lambda)$ (i.e., $\lambda_2\ne\lambda_1$), put
$\varkappa(\lambda_2)=0$.

\begin{theorem}\label{thAsymp1}
Let the lines ${\rm Im\,}\lambda=a_1+1-l-2m$, ${\rm
Im\,}\lambda=a+1-l-2m$ contain no eigenvalues of $\tilde{\cal
L}_1(\lambda)$ and the strip $a_1+1-l-2m<{\rm
Im\,}\lambda<a+1-l-2m$ contain the only eigenvalue $\lambda_{1}$
of $\tilde{\cal L}_1(\lambda)$. Then
 \begin{equation}\label{eqAsymp1}
  \begin{array}{c}
  u(y)=c_{1}u_{1}(y)+c_{2}u_{12}(y)+\hat u(y)\quad (y\in{\cal V}(g_1)\cap
  G).
 \end{array}
 \end{equation}
Here $u_{1}=\{u_{1}^{(k,\zeta)}\}$, where $u_{1}^{(k,\zeta)}$ is
defined by~(\ref{eqPowerSolLambda_n1});
 $u_{12}=\{u_{12}^{(k,\zeta)}\}$, where
$u_{12}^{(k,\zeta)}$ is a linear combination (which will be
strictly defined in the proof below) of the functions
$r^{i\lambda_2}(i\ln r)^q\varphi_{k\zeta q}(\omega)$,
$\varphi_{k\zeta q}\in W_2^{l+2m}(b_1,\ b_2)$, $0\le q\le
k+\varkappa(\lambda_2)$, $(\omega,\ r)$ are polar coordinates with
the pole at $g_1$; $c_{1}=\{c_{1}^{(k,\zeta)}\}$ is a vector of
some constants; $c_2$ is the vector of constants appearing
in~(\ref{eqAsymp2});
 $\hat u\in H_{a_1}^{l+2m}({\cal V}(g_1)\cap G)$\footnotemark{}.
\end{theorem}
\footnotetext{See footnote~\theqountManyLambda~on
page~\pageref{fManyLambda}.}
\begin{proof}
Let $u_{12}=\{u_{12}^{(k,\zeta)}\}$ be a particular solution
(which is defined by Lemma~4.3~\cite{GurPlane}) for the problem
\begin{gather}
 {\cal P}(D_y) u_{12}=0 \quad (y\in K),\label{eqP10Spec}\\
{\cal B}_{1}(D_y)u_{12}=-f_{12},\quad {\cal
B}_{2}(D_y)u_{12}=0,\label{eqB10Spec}
\end{gather}
Here $f_{12}=\{f_{12\mu}^{(k,\zeta)}\}$ is defined
by~(\ref{eqf12}). We remind that each element
$f_{12\mu}^{(k,\zeta)}$ is the linear combination of the functions
$r^{i\lambda_{2}-m_{1\mu}}(i\ln r)^q$, $0\le q\le k$. Therefore,
by Lemma~4.3~\cite{GurPlane}, the particular solution $u_{12}$ has
the form described in the formulation of the theorem. Moreover,
each element $u_{12}^{(k,\zeta)}$ of the vector $u_{12}$ is
uniquely defined if $\lambda_2$ is not an eigenvalue of
$\tilde{\cal L}_1(\lambda)$ (i.e., $\lambda_2\ne\lambda_1$).
Otherwise (i.e., if $\lambda_2=\lambda_1$) it is defined accurate
to an arbitrary linear combination of power
solutions~(\ref{eqPowerSolLambda_n1}) corresponding to the
eigenvalue~$\lambda_2=\lambda_1$. Later on we shall suppose the
particular solution $u_{12}=\{u_{12}^{(k,\zeta)}\}$ being fixed.

Introduce the cut--off function $\eta\in C^\infty({\mathbb R}^2)$
equal to 1 in some neighborhood of the origin and vanishing
outside ${\cal V}(0)$. Put $w=\eta(u-c_2u_{12})$. Since $u$ is a
solution for problem~(\ref{eqP10}), (\ref{eqB10}) and $u_{12}$ is
a solution for problem~(\ref{eqP10Spec}), (\ref{eqB10Spec}), one
can easily check (using Leibnitz's formula) that ${\cal L}_1 w\in
H_{a_1}^{l}(K,\ \gamma)$. Therefore, to conclude the proof, it
remains to apply Theorem~2.2~\cite{GurPlane}, which establishes
the asymptotics of solutions for nonlocal problems in angles.
\end{proof}

\begin{remark}
In fact, the assumption that the line ${\rm Im\,}\lambda=a+1-l-2m$
contains no eigenvalues of $\tilde{\cal L}_1(\lambda)$ is
superfluous. But, using the results of the paper~\cite{GurPlane},
one can show (similarly to Remark~\ref{rSuperfluous2}) that this
assumption does not lead to the loss in generality. Therefore we
remain it since it will be used for studying the adjoint nonlocal
problem and for calculating the coefficients~$c_{1}^{(k,\zeta)}$.
\end{remark}

\begin{remark}\label{rNoEigenval1}
If the strip $a_1+1-l-2m\le{\rm Im\,}\lambda<a+1-l-2m$ has no
eigenvalues of $\tilde{\cal L}_1(\lambda)$, but still has one
(say, $\lambda_2$) of $\tilde{\cal L}_2(\lambda)$,
then~(\ref{eqAsymp1}) will assume the form $
u(y)=c_{2}u_{12}(y)+\hat u(y)$ $(y\in{\cal V}(g_1)\cap G)$. And
only if the mentioned strip has neither eigenvalues of
$\tilde{\cal L}_1(\lambda)$ nor $\tilde{\cal L}_2(\lambda)$, the
solution $u$ will be ``regular'' near the point $g_1$: $u\in
H_{a_1}^{l+2m}({\cal V}(g_1)\cap G)$ for any right--hand side
$\{f,\ f_{\sigma}\}\in H_{a_1}^l(G,\ \Upsilon)$. (Cf.
Remark~\ref{rNoEigenval2}.)
\end{remark}

Theorem~\ref{thAsymp1} shows that the asymptotic behavior of
solutions for problem~(\ref{eqPinG}), (\ref{eqBinG}) near the
point $g_1$ depends on the data of the problem both near the point
$g_1$ itself and near the point $g_2\in G$, which is connected
with $g_1$: $g_2=\Omega_1(g_1)$.

\medskip

{\bf IV.} Quite similarly to the above one can study an
asymptotics of solutions for problem~(\ref{eqPinG}),
(\ref{eqBinG}) near the points $h_\nu$ in terms of the spectral
properties of the operators $\tilde{\cal L}'_\nu(\lambda)$
corresponding to the points $h_\nu$, $\nu=1,\ 2$. The operators
$\tilde{\cal L}'_\nu(\lambda)$ are introduced similarly to the
operators $\tilde{\cal L}_\nu(\lambda)$.

In order not to repeat the analogous computations, we suppose that
the solution $u$ is ``regular'' in some neighborhoods ${\cal
V}(h_\nu)$ of the points $h_\nu$: $u\in H_{a_1}^{l+2m}\big({\cal
V}(h_\nu)\big)$, $\nu=1,\ 2$.

Now we shall formulate the condition that summarize all our
assumptions concerning the spectral properties of the operators
$\tilde{\cal L}_\nu(\lambda)$ and $\tilde{\cal L}'_\nu(\lambda)$.

\begin{condition}\label{condSpectr}
Let the lines ${\rm Im\,}\lambda=a_1+1-l-2m$ and ${\rm
Im\,}\lambda=a+1-l-2m$ contain no eigenvalues of the
operator--valued functions $\tilde{\cal L}_\nu(\lambda)$,
$\tilde{\cal L}_\nu'(\lambda)$; let the strip $a_1+1-l-2m<{\rm
Im\,}\lambda<a+1-l-2m$ contain the only eigenvalue $\lambda_\nu$
of $\tilde{\cal L}_\nu(\lambda)$ and no eigenvalues of
$\tilde{\cal L}_\nu'(\lambda)$, $\nu=1,\ 2$.
\end{condition}

We notice that the assumption concerning the absence of
eigenvalues of $\tilde{\cal L}_\nu'(\lambda)$, $\nu=1,\ 2$, in the
strip $a_1+1-l-2m\le{\rm Im\,}\lambda\le a+1-l-2m$ guarantees the
regularity of solutions in the above sense (see
Remarks~\ref{rNoEigenval2} and~\ref{rNoEigenval1}).

From now on we suppose Condition~\ref{condSpectr} being fulfilled.

\smallskip

In the sequel it will be convenient to have an asymptotics formula
for the solution $u\in H_a^{l+2m}(G)$ to problem~(\ref{eqPinG}),
(\ref{eqBinG}) in the whole domain $G$. To write this formula, we
introduce infinitely smooth functions $\eta_\nu$ with the supports
in ${\cal V}(g_\nu)$ such that $\eta_\nu(y)=1$ in some
neighborhoods of the points $g_\nu$, $\nu=1,\ 2$. Consider the
vector--functions
\begin{equation}\label{eqU12}
 U_{1}=\eta_1 u_{1};\quad
 U_{2}=\eta_2 u_{2} + \eta_1 u_{12}.
\end{equation}
The functions $U_\nu$ are supposed to be defined in the whole
domain $G$, vanishing outside ${\cal V}(g_\nu)$, $\nu=1,\ 2$. Then
Theorems~\ref{thAsymp2} and~\ref{thAsymp1} yield the following
asymptotics of $u\in H_a^{l+2m}(G)$:
\begin{equation}\label{eqAsympG}
  u\equiv \Big(c_1 U_1+c_2 U_2\Big) \Big({\rm mod\,} H_{a_1}^{l+2m}(G)\Big).
\end{equation}

\smallskip

Let us remark for the sequel that the components
$U_{\nu}^{(k,\zeta)}$ of the vector
$U_\nu=\{U_{\nu}^{(k,\zeta)}\}$ are such that
\begin{equation}\label{eqSmoothLU}
{\bf L}U_{\nu}^{(k,\zeta)}\in H_{a_1}^l(G,\ \Upsilon).
\end{equation}

To prove it, we firstly put $\{F,\ F_\sigma\}={\bf
L}U_{1}^{(k,\zeta)}$. Since the support of
$U_{1}^{(k,\zeta)}=\eta_1 u_{1}^{(k,\zeta)}$ is contained in
${\cal V}(g_1)={\cal V}(0)$, we have
\begin{gather}
 {\bf P}(y,\ D_y)\eta_1 u_{1}^{(k,\zeta)}= F(y) \quad (y\in
 K),\notag\\
B_{1}(y,\ D_y)\eta_1 u_{1}^{(k,\zeta)}|_{\gamma_1}= F_1(y) \quad (y\in\gamma_1),\notag\\
B_{2}(y,\ D_y)\eta_1 u_{1}^{(k,\zeta)}|_{\gamma_2}+\big(T_{2}(y,\
D_y)\eta_1 u_{1}^{(k,\zeta)}\big)\big({\cal
G}_{2}y\big)|_{\gamma_2}= F_2(y)
 \quad (y\in\gamma_2).\notag
\end{gather}
But the vector $u_1=\{u_{1}^{(k,\zeta)}\}$
satisfies~(\ref{eqL1u1=0}). Therefore, using Leibnitz's formula,
we obtain $\{F,\ F_\sigma\}\in H_{a_1}^l(G,\ \Upsilon)$.

Now put $\{F,\ F_\sigma\}={\bf L}U_{2}^{(k,\zeta)}$. Similarly to
the above we get
\begin{gather}
{\bf P}(y,\ D_y)\eta_1 u_{12}^{(k,\zeta)} + {\bf P}(y,\ D_y)\eta_2
u_{2}^{(k,\zeta)}= F(y) \quad (y\in K),\notag\\
B_{1}(y,\ D_y)\eta_1 u_{12}^{(k,\zeta)}|_{\gamma_1}+
\big(T_{1}(y,\ D_y)\eta_2
u_{2}^{(k,\zeta)}\big)\big(\Omega_{1}(y)\big)|_{\gamma_1}
= F_1(y) \quad (y\in\gamma_1),\notag\\
B_{2}(y,\ D_y)\eta_1 u_{12}^{(k,\zeta)}|_{\gamma_2}+\big(T_{2}(y,\
D_y)\eta_1 u_{12}^{(k,\zeta)}\big)\big({\cal
G}_{2}y\big)|_{\gamma_2}= F_2(y)
 \quad (y\in\gamma_2).\notag
\end{gather}
But the vector $u_2=\{u_{2}^{(k,\zeta)}\}$
satisfies~(\ref{eqL2u2=0}), and the vector
$u_{12}=\{u_{12}^{(k,\zeta)}\}$ satisfies~(\ref{eqP10Spec}),
(\ref{eqB10Spec}). Therefore, using Leibnitz's formula, we again
obtain $\{F,\ F_\sigma\}\in H_{a_1}^l(G,\ \Upsilon)$.

\section{Index of nonlocal problems}
\phantom{a}~\newline \nopagebreak\label{sectIndex}

{\bf I.} In this section we study some properties of the kernel,
cokernel, and index of the operator ${\bf L}$ corresponding to
nonlocal problem~(\ref{eqPinG}), (\ref{eqBinG}). In particular,
using the asymptotics formula~(\ref{eqAsympG}), we shall obtain a
formula connecting the indices of one and the same
problem~(\ref{eqPinG}), (\ref{eqBinG}), but being considered in
different weighted spaces.

\smallskip

Let $\varkappa_{\nu}$ be a full multiplicity of the eigenvalue
$\lambda_\nu$ of the operator--valued function $\tilde{\cal
L}_\nu(\lambda)$:
$\varkappa_{\nu}=\sum\limits_{\zeta=1}^{J_\nu}\varkappa_{\zeta,\nu}$.
Put $\varkappa=\varkappa_1+\varkappa_2$.

\begin{lemma}\label{lKerL}
Homogeneous problem~(\ref{eqPinG}), (\ref{eqBinG}) can have no
more than $\varkappa$ linearly independent modulo
$H_{a_1}^{l+2m}(G)$ solutions from the space $H_a^{l+2m}(G)$.
\end{lemma}
\begin{proof}
Put the functions $U_{\nu}^{(k,\zeta)}$ ($\nu=1,\ 2$; $\zeta=1,\
\dots,\ J_{\nu}$; $k=0,\ \dots,\ \varkappa_{\zeta,\nu}-1$) in
arbitrary order and denote the elements of the obtained ordered
set by ${\cal U}_1,\ \dots,\ {\cal U}_\varkappa$.

Suppose ${\cal Z}_t\in H_{a_1}^{l+2m}(G)$, $t=1,\ \dots d$, are
linearly independent modulo $H_{a_1}^{l+2m}(G)$ solutions to
homogeneous problem~(\ref{eqPinG}), (\ref{eqBinG}), and
$d>\varkappa$. Then by~(\ref{eqAsympG}) we have
\begin{equation}\label{eqKerL}
 {\cal Z}_t\equiv\Big(\sum\limits_{k=1}^{\varkappa}c_{tk}{\cal U}_k\Big)\ \Big({\rm mod\,}
 H_{a_1}^{l+2m}(G)\Big),\ \quad t=1,\ \dots,\ d,
\end{equation}
where $c_{tk}$ are some constants. Consider the equation for
unknown constants $h_1,\ \dots,\ h_d$:
$$
\sum\limits_{t=1}^{d}h_{t}{\cal Z}_t=0\Big({\rm mod\,}
 H_{a_1}^{l+2m}(G)\Big).
$$
By virtue of~(\ref{eqKerL}) it is equivalent to
$$
\sum\limits_{k=1}^{\varkappa}\Big(\sum\limits_{t=1}^{d}c_{tk}h_t\Big){\cal
U}_k=0 \Big({\rm mod\,}H_{a_1}^{l+2m}(G)\Big).
$$
Since ${\cal U}_1,\ \dots,\ {\cal U}_\varkappa$ are linearly
independent modulo $H_{a_1}^{l+2m}(G)$, the last equation is
equivalent to the system
$$
\sum\limits_{t=1}^{d}c_{tk}h_t=0,\ \quad k=1,\ \dots,\ \varkappa.
$$
By virtue of the inequality $d>\varkappa$, this system necessarily
has a nontrivial solution $(h_1,\ \dots,\ h_d)\ne0$, while we
supposed that ${\cal Z}_1,\ \dots,\ {\cal Z}_d$ were linearly
independent modulo $H_{a_1}^{l+2m}(G)$. This contradiction proves
the lemma.
\end{proof}

Consider the vector ${\cal U}=({\cal U}_1,\ \dots,\ {\cal
U}_\varkappa)^T$. Let ${\cal Z}=({\cal Z}_1,\ \dots,\ {\cal
Z}_d)^T$, $0\le d\le\varkappa$, be a vector, components of which
form a maximal set of solutions to homogeneous
problem~(\ref{eqPinG}), (\ref{eqBinG}) from the space
$H_a^{l+2m}(G)$, linearly independent modulo the space
$H_{a_1}^{l+2m}(G)$ (i.e., a basis modulo $H_{a_1}^{l+2m}(G)$). By
virtue of~(\ref{eqKerL}), we have ${\cal Z}\equiv {\bf C}{\cal U}\
\Big({\rm mod\,} H_{a_1}^{l+2m}(G)\Big)$, where ${\bf C}$ is a
matrix of order $d\times\varkappa$. Rank of ${\bf C}$ equals $d$.
Without loss in generality we assume that ${\bf C}=({\bf C}_1,\
{\bf C}_2)$, where ${\bf C}_1$ is a nonsingular $(d\times
d)$-matrix. Hence ${\bf C}_1^{-1} {\cal Z}\equiv ({\bf I},\ {\bf
C}_1^{-1}{\bf C}_2){\cal U}\ \Big({\rm mod\,}
H_{a_1}^{l+2m}(G)\Big)$, where ${\bf I}$ is the identity $(d\times
d)$-matrix. Therefore we can suppose that
\begin{equation}\label{eqCanonicBasis}
{\cal Z}_t\equiv\Big({\cal U}_t+
\sum\limits_{k=d+1}^{\varkappa}c_{tk}{\cal U}_k\Big)\ \Big({\rm
mod\,} H_{a_1}^{l+2m}(G)\Big),\quad t=1,\ \dots,\ d.
\end{equation}
We shall say that basis~(\ref{eqCanonicBasis}) is canonical. From
now on we fix some canonical basis.

\medskip

{\bf II.} Along with the operator ${\bf L}=\{{\bf P}(y,\ D_y),\
{\bf B}_{\sigma}(y,\ D_y)\}: H_{a_1}^{l+2m}(G)\to  H_{a_1}^l(G,\
\Upsilon)$, we consider the adjoint operator ${\bf
L}^*:H_{a_1}^l(G,\ \Upsilon)^*\to H_{a_1}^{l+2m}(G)^*$ given by
\begin{equation}\label{eqL*inG1}
 <u,\ {\bf L}^*\{v,\ w_\sigma\}>=<{\bf P}(y,\ D_y)u,\ v>+
 \sum\limits_{\sigma=1,2}\sum\limits_{\mu=1}^m
 <{\bf B}_{\sigma\mu}(y,\ D_y)u,\ w_{\sigma\mu}>
\end{equation}
for all $u\in H_{a_1}^{l+2m}(G)$, $\{v,\ w_\sigma\}\in
H_{a_1}^l(G,\ \Upsilon)^*$. Here and below $<\cdot,\ \cdot>$
stands for the sesquilinear form on a pair of corresponding
adjoint spaces.

\begin{lemma}\label{lKerL*}
Let $d$ be a number of elements in basis~(\ref{eqCanonicBasis}).
Then the equation ${\bf L}^*\{v,\ w_{\sigma}\}=0$ has
$\varkappa-d$ solutions from $H_{a_1}^l(G,\ \Upsilon)^*$, linearly
independent modulo $H_{a}^l(G,\ \Upsilon)^*$.
\end{lemma}
\begin{proof}
1) Let $\{\varphi_t,\ \psi_{t,\sigma}\}$, $t=1,\dots,\ q$, be some
basis modulo $H_{a}^l(G,\ \Upsilon)^*$ in the space of solutions
from $H_{a_1}^l(G,\ \Upsilon)^*$ for the equation ${\bf L}^*\{v,\
w_\sigma\}=0$.

Suppose $q<\varkappa-d$. Put ${\cal U}=c_{d+1}{\cal
U}_{d+1}+\dots+c_\varkappa {\cal U}_\varkappa$, where the vector
$(c_{d+1},\ \dots,\ c_\varkappa)$ is a nontrivial solution for the
$q$ linear algebraic equations
\begin{equation}\label{eqKerL*0'}
 <{\bf L}{\cal U},\ \{\varphi_t,\ \psi_{t,\sigma}\}>=0,\quad t=1,\ \dots,\ q
\end{equation}
(notice that, by virtue of~(\ref{eqSmoothLU}), ${\bf L}{\cal U}\in
H_{a_1}^l(G,\ \Upsilon)$ and therefore the forms $<{\bf L}{\cal
U},\ \{\varphi_t,\ \psi_{t,\sigma}\}>$ are well--defined). This
system does have a nontrivial solution since $q<\varkappa-d$.

From~(\ref{eqKerL*0'}) it follows that there exists a solution
$\hat {\cal U}\in H_{a_1}^{l+2m}(G)$ for the equation ${\bf L}\hat
{\cal U}={\bf L}{\cal U}$. Clearly the function ${\cal Z}={\cal
U}-\hat {\cal U}\ne0$ is a solution from $H_{a}^{l+2m}(G)$ for
homogeneous problem~(\ref{eqPinG}), (\ref{eqBinG}), which has the
asymptotics
\begin{equation}\label{eqKerL*0}
 {\cal Z}\equiv\Big(\sum\limits_{k=d+1}^{\varkappa}c_{k}{\cal U}_k\Big)\
\Big({\rm mod\,} H_{a_1}^{l+2m}(G)\Big).
\end{equation}
We claim that the function ${\cal Z}$ is linearly independent  of
 ${\cal Z}_1,\ \dots,\ {\cal Z}_d$, the elements of
basis~(\ref{eqCanonicBasis}) modulo $H_{a_1}^{l+2m}(G)$. Indeed,
suppose that
$$
{\cal Z}\equiv\Big(\sum\limits_{t=1}^{d}h_{t}{\cal Z}_t\Big)\
\Big({\rm mod\,} H_{a_1}^{l+2m}(G)\Big);
$$
then, by virtue of~(\ref{eqCanonicBasis}), we have
$$
{\cal Z}\equiv\Big(\sum\limits_{t=1}^{d}h_{t}{\cal U}_t+
\sum\limits_{k=d+1}^{\varkappa}\bigl(\sum\limits_{t=1}^{d}h_tc_{tk}\bigr){\cal
U}_k\Big)\ \Big({\rm mod\,} H_{a_1}^{l+2m}(G)\Big).
$$
From this, from (\ref{eqKerL*0}), and from the linear independence
of the functions ${\cal U}_1,\ \dots,\ {\cal U}_\varkappa$ modulo
$H_{a_1}^{l+2m}(G)$, it follows that $h_1=\dots=h_d=0$. However,
we assumed ${\cal Z}_1,\ \dots,\ {\cal Z}_d$ were the elements of
basis~(\ref{eqCanonicBasis}) modulo $H_{a_1}^{l+2m}(G)$. This
contradiction proves that $q\ge\varkappa-d$.

2) Suppose $q>\varkappa-d$. Denote by $\{\Phi_h,\
\Psi_{h,\sigma}\}$, $h=1,\ \dots, q$, a system of elements from
$H_{a_1}^l(G,\ \Upsilon)$ being biorthogonal to the system
$\{\varphi_t,\ \psi_{t,\sigma}\}$, $t=1,\ \dots, q$, and
orthogonal to all solutions for the equation ${\bf L}^*\{v,\
w_{\sigma}\}=0$ from $H_{a}^l(G,\ \Upsilon)^*$. Then there exist
solutions $u_h\in H_{a}^{l+2m}(G)$ for the problems ${\bf
L}u_h=\{\Phi_h,\ \Psi_{h,\sigma}\}$, $h=1,\ \dots,\ q$.
Subtracting from $u_h$ (if needed) a linear combination of the
elements ${\cal Z}_1,\ \dots,\ {\cal Z}_d$ forming
basis~(\ref{eqCanonicBasis}), one can make the relations
\begin{equation}\label{eqKerL*1}
 u_h\equiv\Big(\sum\limits_{k=d+1}^{\varkappa}d_{hk}{\cal U}_k\Big)\
\Big({\rm mod\,} H_{a_1}^{l+2m}(G)\Big),\quad h=1,\ \dots,\ q.
\end{equation}
hold.

The functions $u_1,\ \dots,\ u_q$ are linearly independent modulo
$H_{a_1}^{l+2m}(G)$. Indeed, in the opposite case some linear
combination of the functions $u_h$, $h=1,\ \dots,\ q$, would
belong to the space $H_{a_1}^{l+2m}(G)$. Then the corresponding
linear combination of the functions ${\bf L}u_h=\{\Phi_h,\
\Psi_{h,\sigma}\}$, $h=1,\ \dots,\ q$, would be orthogonal to all
the vectors $\{\varphi_t,\ \psi_{t,\sigma}\}$, $t=1,\ \dots, q$.
This would contradict the choice of the functions $\{\Phi_h,\
\Psi_{h,\sigma}\}$, $h=1,\ \dots, q$. From~(\ref{eqKerL*1}) and
from the linear independence of the functions $u_h$, it follows
that $q\le\varkappa-d$. Thus, we necessarily have $q=\varkappa-d$.
\end{proof}

\medskip

{\bf III.} Consider the operators
\begin{align*}
{\bf L}_a &=\{{\bf P}(y,\ D_y),\ {\bf B}_{\sigma}(y,\ D_y)\}: H_a^{l+2m}(G)\to  H_a^l(G,\ \Upsilon),\\
{\bf L}_{a_1} &=\{{\bf P}(y,\ D_y),\ {\bf B}_{\sigma}(y,\ D_y)\}:
H_{a_1}^{l+2m}(G)\to  H_{a_1}^l(G,\ \Upsilon).
 \end{align*}
The operators ${\bf L}_a$ and ${\bf L}_{a_1}$ correspond to one
and the same nonlocal problem~(\ref{eqPinG}), (\ref{eqBinG}), but
they act in the spaces with the different weight constants ($a$
and $a_1$ respectively).

\begin{theorem}\label{thIndL}
The operators ${\bf L}_a$ and ${\bf L}_{a_1}$ are Fredholm, and
the following index formula is valid:
$$
   {\rm ind\,}{\bf L}_a={\rm ind\,}{\bf L}_{a_1}+\varkappa.
$$
\end{theorem}
\begin{proof}
By Theorem~3.4~\cite{SkMs86}, the operators  ${\bf L}_a$ and ${\bf
L}_{a_1}$ are Fredholm\footnote{More precisely, we use the
generalization of Theorem~3.4~\cite{SkMs86} for the case of
transformations $\Omega_1,\ \Omega_2$ consisting near $g_1$ and
$h_1$ not only of a rotation but of an expansion, too; see
also~\cite{SkDu91}.}. By Lemma~\ref{lKerL}, we have ${\rm
dim\,}{\rm ker\,}{\bf L}_a={\rm dim\,}{\rm ker\,}{\bf L}_{a_1}+d$.
Then by Lemma~\ref{lKerL*}, we have ${\rm dim\,}{\rm ker\,}{\bf
L}^*_a={\rm dim\,}{\rm ker\,}{\bf L}^*_{a_1}-(\varkappa-d)$. Hence
${\rm ind\,}{\bf L}_a={\rm dim\,}{\rm ker\,}{\bf L}_a-{\rm
dim\,}{\rm ker\,}{\bf L}^*_a={\rm dim\,}{\rm ker\,}{\bf
L}_{a_1}-{\rm dim\,}{\rm ker\,}{\bf L}^*_{a_1}+\varkappa= {\rm
ind\,}{\bf L}_{a_1}+\varkappa$.
\end{proof}

\begin{remark}
Theorem~\ref{thIndL} remains true without the assumption
$a-a_1<1$, too. Indeed, one can always choose numbers
$a=a^0>a^1>\dots>a^M=a_1$ such that $0<a^{i}-a^{i+1}<1$ and the
lines ${\rm Im\,}\lambda=a^i+1-l-2m$ do not contain eigenvalues of
$\tilde{\cal L}_\nu(\lambda)$, $\nu=1,\ 2$. Applying
Theorem~\ref{thIndL} subsequently to the pairs of the operators
 \begin{align*}
{\bf L}_{a_i}&=\{{\bf P}(y,\ D_y),\ {\bf B}_{\sigma}(y,\ D_y)\}:
H_{a_i}^{l+2m}(G)\to  H_{a_i}^l(G,\ \Upsilon),\\
{\bf L}_{a_{i+1}}&=\{{\bf P}(y,\ D_y),\ {\bf B}_{\sigma}(y,\
D_y)\}: H_{a_{i+1}}^{l+2m}(G)\to  H_{a_{i+1}}^l(G,\ \Upsilon)
 \end{align*}
we get the formula ${\rm ind\,}{\bf L}_a={\rm ind\,}{\bf
L}_{a_1}+\varkappa$, where $\varkappa$ is the sum of full
multiplicities of all eigenvalues of $\tilde{\cal L}_1(\lambda)$
and $\tilde{\cal L}_2(\lambda)$ contained in the strip
$a_1+1-l-2m<{\rm Im\,}\lambda<a+1-l-2m$.
\end{remark}

\section{Asymptotics of solutions for adjoint nonlocal problems}
\phantom{a}~\newline \nopagebreak\label{sectAsymp*}

{\bf I.} In this section we shall obtain an asymptotics near the
set ${\cal K}$ for solutions to the problem, adjoint
to~(\ref{eqPinG}), (\ref{eqBinG}). The results of this section
will be applied to calculating the coefficients $c_{\nu}^{(k,j)}$
in~(\ref{eqAsympG}).

\smallskip

Notice that the approach to the study of adjoint nonlocal problems
has been suggested by the author in~\cite{GurDAN, GurGiess,
GurPlane}. In the papers~\cite{GurDAN, GurGiess}, the solvability
and smoothness of solutions for model adjoint nonlocal problems in
plane and dihedral angles were studied. The paper~\cite{GurPlane}
deals with an asymptotics of solutions for model nonlocal problems
in plane angles and in~${\mathbb R}^2\setminus\{0\}$. In the
present work we essentially use both the ideology of the
papers~\cite{GurDAN, GurGiess} and the results of the
paper~\cite{GurPlane}.

\smallskip

Since we suppose Condition~\ref{condSpectr} being fulfilled, it
suffices to obtain appropriate asymptotics formulas only near the
points $g_1$ and $g_2$.

Along with formula~(\ref{eqL*inG1}) we will use another one for
the definition of the adjoint operator. To write this formula, we
introduce the following denotation. For any smooth curve
$\Upsilon\subset\bar G$ and any distribution $w\in
H_{a_1}^{k-1/2}(\Upsilon)^*$ we denote by $w\cdot\delta_\Upsilon$
the distribution from $H_{a_1}^k(G)^*$ given by
\begin{equation}\label{eqwdelta}
 <u,\ w\cdot\delta_\Upsilon>_G=<u|_{\Upsilon},\ w>_\Upsilon\
 \mbox{for all } u\in H_{a_1}^k(G)\footnotemark.
\end{equation}
\footnotetext{In this section, for clearness, we denote
sesquilinear forms on the pairs of adjoint spaces $H_{a_1}^k(G)$,
$H_{a_1}^k(G)^*$ and $H_{a_1}^{k-1/2}(\Upsilon)$,
$H_{a_1}^{k-1/2}(\Upsilon)^*$ by $<\cdot,\ \cdot>_G$ and $<\cdot,\
\cdot>_\Upsilon$ respectively.} Clearly the support of the
distribution $w\cdot\delta_\Upsilon$ is contained in
$\bar\Upsilon$. Similarly one can define a distribution
$w\cdot\delta_{\gamma}\in H_{a_1}^k(K)^*$, where
$\gamma=\{y\in{\mathbb R}^2:\ r>0,\ \omega=b\}$ ($b_1\le b\le
b_2$) and $w\in H_{a_1}^{k-1/2}(\gamma)^*$.

Denote by ${\bf P}^*(y,\ D_y)$, $B_{\sigma}^*(y,\
D_y)=\{B_{\sigma\mu}^*(y,\ D_y)\}_{\mu=1}^m$, $T_{\sigma}^*(y,\
D_y)=\{T_{\sigma\mu}^*(y,\ D_y)\}_{\mu=1}^m$ the operators,
formally adjoint to ${\bf P}(y,\ D_y)$, $B_{\sigma}(y,\
D_y)=\{B_{\sigma\mu}(y,\ D_y)\}_{\mu=1}^m$, $T_{\sigma}(y,\
D_y)=\{T_{\sigma\mu}(y,\ D_y)\}_{\mu=1}^m$ respectively.

For any distribution $w_{\sigma\mu}\in
H_{a_1}^{l+2m-m_{\sigma\mu}-1/2}(\Upsilon_{\sigma})$, we consider
the distribution $w_{\sigma\mu}^\Omega\in
H_{a_1}^{l+2m-m_{\sigma\mu}-1/2}\big(\Omega_\sigma(\Upsilon_{\sigma})\big)^*$
 given by
\begin{equation}\label{eqwOmega}
<\psi,\ w_{\sigma\mu}^\Omega>_{\Omega_\sigma(\Upsilon_{\sigma})}
=<\psi\big(\Omega_{\sigma}(\cdot)\big),\
w_{\sigma\mu}>_{\Upsilon_\sigma}
\end{equation}
for all $\psi\in
H_{a_1}^{l+2m-m_{\sigma\mu}-1/2}\big(\Omega_\sigma(\Upsilon_{\sigma})\big)$.

We claim that the adjoint operator ${\bf L}^*:H_{a_1}^l(G,\
\Upsilon)^*\to H_{a_1}^{l+2m}(G)^*$ can be defined by the formula
\begin{equation}\label{eqL*inG2}
{\bf L}^*\{v,\ w_\sigma\}={\bf P}^*(y,\
D_y)v+\sum\limits_{\sigma=1,2} B_{\sigma}^*(y,\
D_y)(w_{\sigma}\cdot\delta_{\Upsilon_\sigma}) +T_{\sigma}^*(y,\
D_y)(w_{\sigma}^\Omega\cdot\delta_{\Omega_\sigma(\Upsilon_\sigma)}).
\end{equation}
Here and further $w_{\sigma}=\{w_{\sigma\mu}\}_{\mu=1}^m$,
$w_{\sigma}^\Omega=\{w_{\sigma\mu}^\Omega\}_{\mu=1}^m$,
\begin{gather}
B_{\sigma}^*(y,\ D_y)
(w_{\sigma}\cdot\delta_{\Upsilon_\sigma})=\sum\limits_{\mu=1}^m
B_{\sigma\mu}^*(y,\
D_y)(w_{\sigma\mu}\cdot\delta_{\Upsilon_\sigma}),\notag\\
T_{\sigma}^*(y,\ D_y)
(w_{\sigma}^\Omega\cdot\delta_{\Omega_\sigma(\Upsilon_\sigma)})=
\sum\limits_{\mu=1}^m T_{\sigma\mu}^*(y,\
D_y)(w_{\sigma\mu}^\Omega\cdot\delta_{\Omega_\sigma(\Upsilon_\sigma)}).\notag
\end{gather}

Indeed, using definition~(\ref{eqL*inG1}) of the adjoint operator
${\bf L}^*$ and then relations~(\ref{eqwOmega})
and~(\ref{eqwdelta}), we get (omitting $(y,\ D_y)$ for short)
\begin{multline}\notag
 <u,\ {\bf L}^*\{v,\ w_{\sigma}\}>=<{\bf P}u,\ v>_G+
\sum\limits_{\sigma=1,2}\sum\limits_{\mu=1}^m\Big(
<B_{\sigma\mu}u|_{\Upsilon_{\sigma}},\
w_{\sigma\mu}>_{\Upsilon_{\sigma}} +\\
<\big(T_{\sigma\mu}u\big)\big(\Omega_\sigma(\cdot)\big)|_{\Upsilon_{\sigma}},\
w_{\sigma\mu}>_{\Upsilon_{\sigma}}\Big) =<{\bf P}u,\ v>_G+\\
\sum\limits_{\sigma=1,2}\sum\limits_{\mu=1}^m\Big(
<B_{\sigma\mu}u,\ w_{\sigma\mu}\cdot\delta_{\Upsilon_\sigma}>_{G}
+ <T_{\sigma\mu}u,\
w_{\sigma\mu}^\Omega\cdot\delta_{\Omega_\sigma(\Upsilon_\sigma)}>_{G}\Big)
\end{multline}
for all $u\in H_{a_1}^{l+2m}(G)$, which yields~(\ref{eqL*inG2}).

We are to study an asymptotics of a given solution $\{v,\
w_{\sigma}\}\in H_{a_1}^l(G,\ \Upsilon)^*$ for the problem
\begin{equation}\label{eqL*G}
 {\bf L}^*\{v,\ w_{\sigma}\}=\Psi,
\end{equation}
supposing $\Psi\in H_a^{l+2m}(G)^*$.

For this purpose, parallel to the operator ${\bf L}^*$, we
consider the auxiliary operator
\begin{multline}\notag
{\bf L}_\Omega^*:
H_{a_1}^l(G)^*\times\prod\limits_{\sigma=1,2}\prod\limits_{\mu=1}^m
\Big(H_{a_1}^{l+2m-m_{\sigma\mu}-1/2}(\Upsilon_{\sigma})^*\times\\
\times
H_{a_1}^{l+2m-m_{\sigma\mu}-1/2}\big(\Omega_\sigma(\Upsilon_{\sigma})\big)^*\Big)
\to H_{a_1}^{l+2m}(G)^*
\end{multline}
given by
\begin{multline}\label{eqL*inGAuxil}
{\bf L}_\Omega^*\{v,\ w_\sigma,\ w_\sigma'\}={\bf P}^*(y,\
D_y)v+\\
\sum\limits_{\sigma=1,2}\Big( B_{\sigma}^*(y,\
D_y)(w_{\sigma}\cdot\delta_{\Upsilon_\sigma}) +T_{\sigma}^*(y,\
D_y)(w_{\sigma}'\cdot\delta_{\Omega_\sigma(\Upsilon_\sigma)})\Big).
\end{multline}

Such an auxiliary operator was used in the papers~\cite{GurMatZam,
GurGiess} for the study of model nonlocal problems in angles. It
also turns out to be very useful in our case. On the one hand, the
operator ${\bf L}_\Omega^*$ is not a nonlocal one since the
functions $w_{\sigma}$ and $w_{\sigma}'$ are not connected each
with another by nonlocal transformations $\Omega_\sigma$
in~(\ref{eqL*inGAuxil}). This will allow to use Leibnitz's formula
when necessary. On the other hand, solutions for the problems
corresponding to ${\bf L}^*$ and ${\bf L}_\Omega^*$ are related in
the following way. If $\{v,\ w_\sigma\}$ is a solution to
problem~(\ref{eqL*G}), then the distribution $\{v,\ w_\sigma,\
w_\sigma^\Omega\}$ is a solution to the problem
\begin{equation}\label{eqLOmega*G}
 {\bf L}^*_\Omega\{v,\ w_\sigma,\ w_\sigma^\Omega\}=\Psi,
\end{equation}
where $w_\sigma^\Omega=\{w_{\sigma\mu}^\Omega\}_{\mu=1}^m$ is
defined by~(\ref{eqwOmega}).

So, investigating the asymptotics of the solution $\{v,\
w_{\sigma},\ w_\sigma^\Omega\}$ for problem~(\ref{eqLOmega*G}) is
equivalent to investigating the asymptotics of the solution $\{v,\
w_{\sigma}\}$ for problem~(\ref{eqL*G}).

Our plan is this. At first we will multiply the distribution
$\{v,\ w_\sigma,\ w_\sigma^\Omega\}$ by the cut--off function
$\eta_1$ and get a corresponding problem near $g_1$. Since the
operator ${\bf L}_\Omega^*$ is not a nonlocal one, applying
Leibnitz's formula we will show that ${\bf L}_\Omega^*\eta_1\{v,\
w_{\sigma},\ w_\sigma^\Omega\}\in H_a^{l+2m}(G)^*$. Therefore we
will arrive at the model adjoint problem in the angle $K$ with a
``regular'' right--hand side. Using the results
of~\cite{GurPlane}, we will obtain the asymptotics near $g_1$.
Then we will multiply the distribution $\{v,\ w_\sigma,\
w_\sigma^\Omega\}$ by $\eta_2$ and get a corresponding problem
near $g_2$. We will arrive at the model adjoint problem in
${\mathbb R}^2$. But in this case a right--hand side will be a sum
of ``regular'' and ``special'' distributions. The asymptotics of
the ``special'' one will be defined by the asymptotics of $w_{1}$
near $g_1$, which will have been known from the first step.
Further application of the results of~\cite{GurPlane} will allow
to get the asymptotics near $g_2$.

\smallskip

Thus let us multiply $\{v,\ w_\sigma,\ w_\sigma^\Omega\}$ by
$\eta_1$. Notice that ${\rm
supp\,}\eta_1\cap\overline{\Omega_1(\Upsilon_1)}=\varnothing$ (see
Fig.~\ref{figDomG}) and ${\rm supp\,}
(w_{1}\cdot\delta_{\Omega_1(\Upsilon_1)})\subset\overline{\Omega_1(\Upsilon_1)}$.
Therefore
$\eta_1w_{1}^\Omega\cdot\delta_{\Omega_1(\Upsilon_1)}=0$. From
this and from~(\ref{eqL*inGAuxil}), it follows that
\begin{multline}\label{eqLOmega*1}
{\bf L}_\Omega^*\eta_1\{v,\ w_\sigma,\ w_\sigma^\Omega\}={\bf
P}^*(y,\ D_y)\eta_1v + B_{1}^*(y,\
D_y)(\eta_1w_{1}\cdot\delta_{\Upsilon_1})+\\
B_{2}^*(y,\ D_y)(\eta_1w_{2}\cdot\delta_{\Upsilon_2}) +T_{2}^*(y,\
D_y)(\eta_1w_{2}^\Omega\cdot\delta_{\Omega_2(\Upsilon_2)}).
\end{multline}

Let us show that the distribution $\eta_1\{v,\ w_\sigma\}$
satisfies the model adjoint problem in the angle $K$
\begin{equation}\label{eqL*K}
 {\cal L}_1^*\eta_1\{v,\ w_\sigma\}=\hat\Psi,
\end{equation}
where ${\cal L}_1^*: H_{a_1}^{l}(K,\ \gamma)^*\to
H_{a_1}^{l+2m}(K)^*$ is the operator, adjoint to ${\cal L}_1:
H_{a_1}^{l+2m}(K)\to H_{a_1}^{l}(K,\ \gamma)$; $\hat\Psi\in
H_a^{l+2m}(K)^*$.

From~(\ref{eqLOmega*G}) and Leibnitz's formula, it follows that,
on the one hand,
\begin{equation}\label{eqeta1PsihatPsi1}
{\bf L}_\Omega^*\eta_1\{v,\ w_\sigma,\
w_\sigma^\Omega\}=\eta_1\Psi+\hat\Psi_1,
\end{equation}
where $\hat\Psi_1\in H_a^{l+2m}(G)^*$ and ${\rm
supp\,}\hat\Psi_1\subset{\cal V}(g_1)$. On the other hand, the
function $\eta_1(y)-\eta_1\big(\Omega_2^{-1}(y)\big)$ is equal to
0 near $g_1$ and has a support inside ${\cal V}(g_1)$ (here we are
to suppose ${\rm supp\,}\eta_1$ is so small that ${\rm
supp\,}\eta_1\big(\Omega_2^{-1}(\cdot)\big)\subset{\cal V}(g_1)$).
Hence,
$$
T_{2}^*(y,\ D_y)
(\eta_1w_{2}^\Omega\cdot\delta_{\Omega_2(\Upsilon_2)})-
T_{2}^*(y,\ D_y)
\big(\eta_1(\Omega_2^{-1}(\cdot))w_{2}^\Omega\cdot
\delta_{\Omega_2(\Upsilon_2)}\big)\in H_a^{l+2m}(G)^*
$$
and has a support inside ${\cal V}(g_1)$. This,
(\ref{eqLOmega*1}), and~(\ref{eqeta1PsihatPsi1}) imply
 \begin{multline}\notag
{\bf P}^*(y,\ D_y)\eta_1v + B_{1}^*(y,\
D_y)(\eta_1w_{1}\cdot\delta_{\Upsilon_1})+\\
B_{2}^*(y,\ D_y)(\eta_1w_{2}\cdot\delta_{\Upsilon_2}) +T_{2}^*(y,\
D_y)
(\eta_1(\Omega_2^{-1}(\cdot))w_{2}^\Omega\cdot\delta_{\Omega_2(\Upsilon_2)})=
\hat\Psi_2,
\end{multline}
where $\hat\Psi_2\in H_a^{l+2m}(G)^*$ and ${\rm
supp\,}\hat\Psi_2\subset{\cal V}(g_1)$.

Let ${\cal P}^*(D_y),\ B_{2}^*(D_y),\ T_{2}^*(D_y)$ be the
principal homogeneous parts of the operators ${\bf P}^*(g_1,\
D_y),\ B_{2}^*(g_1,\ D_y),\ T_{2}^*(g_1,\ D_y)$ respectively.
Then, using Leibnitz's formula, we finally get
\begin{multline}\label{eqL*K'}
{\cal P}^*(D_y)\eta_1v + B_{1}^*(D_y)(\eta_1w_{1}\cdot\delta_{\gamma_1})+\\
B_{2}^*(D_y)(\eta_1w_{2}\cdot\delta_{\gamma_2}) +T_{2}^*(D_y)
(\eta_1({\cal G}_2^{-1}\cdot)w_{2}^\Omega\cdot\delta_{{\cal
G}_2(\gamma_2)})= \hat\Psi,
\end{multline}
where $\hat\Psi\in H_a^{l+2m}(K)^*$ and ${\rm
supp\,}\hat\Psi\subset{\cal V}(0)$. Here we also took into account
that near the point $g_1=0$ the domain $G$ and the curves
$\Upsilon_\sigma$ coincide with the angle $K$ and the arms
$\gamma_\sigma$ respectively, while the transformation $\Omega_2$
coincides with the linear operator ${\cal G}_2$. But it is  easily
seen that equality~(\ref{eqL*K'}) is quite the same as
equality~(\ref{eqL*K}). Indeed, the only not evident identity one
should check is
 $$
 <u,\ T_{2}^*(D_y) (\eta_1({\cal
G}_2^{-1}\cdot)w_{2}^\Omega\cdot\delta_{{\cal
G}_2(\gamma_2)})>_K=<(T_{2}(D_y)u)({\cal G}_2\cdot)|_{\gamma_2},\
\eta_1w_2>_{\gamma_2},
$$ which follows from
 \begin{multline}\notag
 <u,\ T_{2}^*(D_y) \big(\eta_1({\cal
G}_2^{-1}\cdot)w_{2}^\Omega\cdot\delta_{{\cal
G}_2(\gamma_2)}\big)>_K= \\
<\eta_1({\cal G}_2^{-1}\cdot)T_{2}(D_y)u|_{{\cal G}_2(\gamma_2)},\
w_2^\Omega>_{{\cal G}_2(\gamma_2)} =<\eta_1(T_{2}(D_y)u)({\cal
G}_2\cdot)|_{\gamma_2},\ w_2>_{\gamma_2}.
\end{multline}
Here we subsequently used~(\ref{eqwdelta}) and~(\ref{eqwOmega}).

\smallskip

Applying the results of the paper~\cite{GurPlane} to
equality~(\ref{eqL*K}), we shall now obtain the asymptotics of the
distribution $\eta_1\{v,\ w_\sigma\}$. Introduce some denotation.
Put
\begin{equation}\label{eqPowerSol1*}
\begin{array}{c}
v_1^{(k,\zeta)}=
r^{i\bar\lambda_{1}+2m-2}\sum\limits_{q=0}^k\frac{\displaystyle
1}{\displaystyle q!}(i\ln r)^q\psi_{1}^{(k-q,\zeta)},\\
w_{1,\sigma\mu}^{(k,\zeta)}=
r^{i\bar\lambda_{1}+m_{\sigma\mu}-1}\sum\limits_{q=0}^k\frac{\displaystyle
1}{\displaystyle q!}(i\ln r)^q\chi_{1,\sigma\mu}^{(k-q,\zeta)}.
\end{array}
\end{equation}
Here
$$
\Big\{\{\psi_{1}^{(0,\zeta)},\ \chi_{1,\sigma\mu}^{(0,\zeta)}\},\
\dots,\ \{\psi_{1}^{(\varkappa_{\zeta,1}-1,\zeta)},\
\chi_{1,\sigma\mu}^{(\varkappa_{\zeta,1}-1,\zeta)}\}:\
 \zeta=1,\ \dots,\ J_{1} \Big\}
$$
are Jordan chains of the operator $\tilde{\cal L}_1^*(\lambda)$
(adjoint to $\tilde{\cal L}_1(\bar\lambda)$) corresponding to the
eigenvalue $\bar\lambda_{1}$ and forming a canonical system. These
chains are supposed (see Lemma~3.2~\cite{GurPlane}) to satisfy the
following condition of biorthogonality and normalization with
respect to the Jordan chains~(\ref{eqJordan1}):
\begin{multline}\label{eqNorm1}
\sum\limits_{p=0}^\nu\sum\limits_{q=0}^k\frac{\displaystyle
1}{\displaystyle (\nu+k+1-p-q)!}
<\partial_\lambda^{\nu+k+1-p-q}\tilde{\cal
L}_1(\lambda_1)\varphi_1^{(q,\xi)},\\
 \{\psi_1^{(p,\zeta)},\
\chi_{1,\sigma\mu}^{(p,\zeta)}\}>
 =
\delta_{\xi,\zeta}\delta_{\varkappa_{\xi,1}-k-1,\nu}.
\end{multline}
Here $\zeta,\ \xi=1,\ \dots,\ J_1$; $\nu=0,\ \dots,\
\varkappa_{\zeta,1}-1$; $k=0,\ \dots,\ \varkappa_{\xi,1}-1$;
$\delta_{\xi,\zeta}$ is the Kronecker symbol.

Analogously to section~\ref{sectAsymp}, we introduce the vectors
$w_{1,\sigma}^{(k,\zeta)}=\{w_{1,\sigma\mu}^{(k,\zeta)}\}_{\mu=1}^m$
and $\{v_1,\ w_{1,\sigma}\}=\{v_{1}^{(k,\zeta)},\
w_{1,\sigma}^{(k,\zeta)}\}$.

We remark that, by Lemma~3.1~\cite{GurPlane}, the distributions
$\{v_{1}^{(k,\zeta)},\ w_{1,\sigma}^{(k,\zeta)}\}$ satisfy the
homogeneous equation ${\cal L}_1^*\{v_{1}^{(k,\zeta)},\
w_{1,\sigma}^{(k,\zeta)}\}=0$.

Now from equality~(\ref{eqL*K}) and Theorem~4.2~\cite{GurPlane} we
get the following result.
\begin{theorem}\label{thAsymp*1}
Let $\{v,\ w_{\sigma}\}\in H_{a_1}^l(G,\ \Upsilon)^*$ be a
solution for equation~(\ref{eqL*G}) with a right--hand side
$\Psi\in H_a^{l+2m}(G)^*$. Then the following asymptotics formula
is valid:
 \begin{equation}\label{eqAsymp*1}
   \eta_1\{v,\ w_{\sigma}\}\equiv
  d_{1}\eta_1\{v_{1},\ w_{1,\sigma}\}\ \Big({\rm mod\,} H_a^{l}(G,\ \Upsilon)^*\Big),
 \end{equation}
where $\{v_1,\ w_{1,\sigma}\}=\{v_{1}^{(k,\zeta)},\
w_{1,\sigma}^{(k,\zeta)}\}$ is defined by~(\ref{eqPowerSol1*}),
$d_{1}=\{d_{1}^{(k,\zeta)}\}$ is a vector of some constants.
\end{theorem}

\medskip

{\bf II.} Now let us study the asymptotics of the solution $\{v,\
w_\sigma\}$ for equation~(\ref{eqL*G}) near the point $g_2$. As we
mentioned above, in this case we will arrive at the model adjoint
problem in ${\mathbb R}^2\setminus\{g_2\}$. A right--hand side of
the equation obtained will be a sum of ``regular'' and ``special''
distributions. The asymptotics of the latter one will be defined
by the asymptotics of $w_{\sigma}$ near $g_1$, which is already
known (see Theorem~\ref{thAsymp*1}).

\smallskip

We multiply $\{v,\ w_\sigma,\ w_\sigma^\Omega\}$ by $\eta_2$.
Since the support of $\eta_2$ is contained in ${\cal V}(g_2)$ and
therefore does not intersect with $\bar\Upsilon_1$,
$\bar\Upsilon_2$, and $\overline{\Omega_2(\Upsilon_2)}$ (see
Fig.~\ref{figDomG}), we have
$\eta_2w_{1}\cdot\delta_{\Upsilon_1}=0$,
$\eta_2w_{2}\cdot\delta_{\Upsilon_2}=0$, and
$\eta_2w_{2}^\Omega\cdot\delta_{\Omega_2(\Upsilon_2)}=0$.
Combining this with~(\ref{eqL*inGAuxil}), we get
\begin{equation}\label{eqLOmega*2}
{\bf L}_\Omega^*\eta_2\{v,\ w_\sigma,\ w_\sigma^\Omega\}={\bf
P}^*(y,\ D_y)\eta_2v + T_{1}^*(y,\
D_y)(\eta_2w_{1}^\Omega\cdot\delta_{\Omega_1(\Upsilon_1)}).
\end{equation}
From~(\ref{eqLOmega*G}) and Leibnitz's formula, it follows that $
{\bf L}_\Omega^*\eta_2\{v,\ w_\sigma,\
w_\sigma^\Omega\}=\eta_2\Psi+\hat\Psi_2, $ where $\hat\Psi_2\in
H_a^{l+2m}(G)^*$ and ${\rm supp\,}\hat\Psi_2\subset{\cal V}(g_2)$.

Let ${\cal P}^*(D_y)$, $T_{1}^*(D_y)$ be the principal homogeneous
parts of the operators ${\bf P}^*(g_2,\ D_y)$, $T_{1}^*(g_2,\
D_y)$ respectively. Then, analogously to the above, we derive that
\begin{equation}\label{eqL*2}
 \begin{array}{c}
 {\cal P}^*(D_y)\eta_2 v
 =-T_{1}^*(D_y)(\eta_2w_{1}^\Omega\cdot\delta_{\Omega_1(\Upsilon_1)})
 +\hat\Psi_3,
 \end{array}
\end{equation}
where $\hat\Psi_3\in H_a^{l+2m}(G)^*$ and ${\rm
supp\,}\hat\Psi_3\subset{\cal V}(g_2)$.

Let ${\cal L}_2^*: H_{a_1}^{l}({\mathbb R}^2)^*\to
H_{a_1}^{l+2m}({\mathbb R}^2)^*$ be the operator, adjoint to
${\cal L}_2: H_{a_1}^{l+2m}({\mathbb R}^2)\to H_{a_1}^{l}({\mathbb
R}^2)$.

From definition~(\ref{eqwOmega}) of the distribution
$w_{1\mu}^\Omega$ and from asymptotics formula~(\ref{eqAsymp*1}),
it follows that $\eta_2w_{1\mu}^\Omega$ is a linear combination of
the functions $r^{i\bar\lambda_{1}+m_{1\mu}-1}(i\ln r)^q$ modulo
$H_a^{l+2m-m_{1\mu}-1/2}(\Upsilon_1)$, where $r$ is a polar radius
of polar coordinates with the pole at $g_2$. Since
$T_{1\mu}^*(D_y)$ is a homogeneous operator of order~$m_{1\mu}$,
we can write~(\ref{eqL*2}) (taking into account~(\ref{eqAsymp*1}))
in the form
\begin{equation}\label{eqL*20}
 \begin{array}{c}
 {\cal L}_2^*\eta_2 v
 =-\eta_2d_1\Psi_{21}+\hat\Psi.
 \end{array}
\end{equation}
Here $\hat\Psi\in H_a^{l+2m}({\mathbb R}^2)^*$ and ${\rm
supp\,}\hat\Psi\subset{\cal V}(g_2)$;
$\Psi_{21}=\{\Psi_{21}^{(k,\zeta)}\}$, where
$\Psi_{21}^{(k,\zeta)}$ is a linear combination of the
distributions $r^{i\bar\lambda_{1}-2}(i\ln
r)^q\Psi_{21q}^{(k,\zeta)}$, $0\le q\le k$,
$\Psi_{21q}^{(k,\zeta)}\in W_{2,2\pi}^{l+2m}(0,\ 2\pi)^*$, $r$ is
a polar radius of polar coordinates with the pole at
$g_2$\footnote{The distribution $\eta_2r^{i\bar\lambda_{1}-2}(i\ln
r)^q\Psi_{21q}^{(k,\zeta)}\in H_{a_1}^{l+2m}({\mathbb R}^2)^*$ is
given by
$$<u,\
r^{i\bar\lambda_{1}-2}(i\ln r)^q\Psi_{21q}^{(k,\zeta)}>=
\int\limits_0^\infty<\eta_2(\cdot,\ r)u(\cdot,\ r),\
\Psi_{21q}^{(k,\zeta)}>_{(0,\
2\pi)}\overline{r^{i\bar\lambda_{1}-2} (i\ln r)^q}\,r\,dr\
$$
for all $u\in H_{a_1}^{l+2m}({\mathbb R}^2)$.}; $d_1$ is the
vector of constants from~(\ref{eqAsymp*1}).

Thus we see that~(\ref{eqL*20}) is a model adjoint problem in
${\mathbb R}^2$ with the right--hand side being the sum of the
``regular'' distribution $\hat\Psi$ and the ``special''
distribution $\Psi_{21}$. The asymptotics of $\Psi_{12}$ is
defined by the asymptotics of the solution $w_1$ near the point
$g_1$, i.e., by the functions $w_{1,\sigma\mu}^{(k,\zeta)}$
(see~(\ref{eqPowerSol1*})).

\smallskip

Applying the results of the paper~\cite{GurPlane} to
equality~(\ref{eqL*20}), we shall now obtain the asymptotics of
the distribution $\eta_2 v$. Introduce some denotation. Put
\begin{equation}\label{eqPowerSol2*}
v_2^{(k,\zeta)}=
r^{i\bar\lambda_{2}+2m-2}\sum\limits_{q=0}^k\frac{\displaystyle
1}{\displaystyle q!}(i\ln r)^q\psi_{2}^{(k-q,\zeta)}.
\end{equation}
Here
$$
\Big\{\psi_{2}^{(0,\zeta)},\ \dots,\
\psi_{2}^{(\varkappa_{\zeta,2}-1,\zeta)}:\
 \zeta=1,\ \dots,\ J_{2} \Big\}
$$
are Jordan chains of the operator $\tilde{\cal L}_2^*(\lambda)$
(adjoint to $\tilde{\cal L}_2(\bar\lambda)$) corresponding to the
eigenvalue $\bar\lambda_{2}$ and forming a canonical system. These
chains are supposed (see~\cite{GurPlane}) to satisfy the following
condition of biorthogonality and normalization with respect to the
Jordan chains~(\ref{eqJordan2}):
\begin{multline}\label{eqNorm2}
\sum\limits_{p=0}^\nu\sum\limits_{q=0}^k\frac{\displaystyle
1}{\displaystyle (\nu+k+1-p-q)!}
<\partial_\lambda^{\nu+k+1-p-q}\tilde{\cal
L}_2(\lambda_2)\varphi_2^{(q,\xi)},\\
 \psi_2^{(p,\zeta)}>
 =
\delta_{\xi,\zeta}\delta_{\varkappa_{\xi,2}-k-1,\nu}.
\end{multline}
Here $\zeta,\ \xi=1,\ \dots,\ J_2$; $\nu=0,\ \dots,\
\varkappa_{\zeta,2}-1$; $k=0,\ \dots,\ \varkappa_{\xi,2}-1$.

Analogously to section~\ref{sectAsymp}, we introduce the vector
$v_2=\{v_{2}^{(k,\zeta)}\}$.

We remark that, according to~\cite[section~5]{GurPlane}, the
distributions $v_{2}^{(k,\zeta)}$ satisfy the homogeneous equation
${\cal L}_2^*v_{2}^{(k,\zeta)}=0$.

If $\bar\lambda_1$ is an eigenvalue of $\tilde{\cal
L}^*_2(\lambda)$ (i.e., $\bar\lambda_1=\bar\lambda_2$), then
denote by $\varkappa(\bar\lambda_1)$ the greatest of partial
multiplicities of $\bar\lambda_1$. If $\bar\lambda_1$ is not an
eigenvalue of $\tilde{\cal L}^*_2(\lambda)$ (i.e.,
$\bar\lambda_1\ne\bar\lambda_2$), put
$\varkappa(\bar\lambda_1)=0$.

\begin{theorem}\label{thAsymp*2}
Let $\{v,\ w_{\sigma}\}\in H_{a_1}^l(G,\ \Upsilon)^*$ be a
solution for equation~(\ref{eqL*G}) with a right--hand side
$\Psi\in H_a^{l+2m}(G)^*$. Then the following asymptotics formula
is valid:
 \begin{equation}\label{eqAsymp*2}
   \eta_2 v\equiv
  \Big(d_{2}\eta_2v_{2}+d_{1}\eta_2v_{21}\Big) \Big({\rm mod\,} H_a^{l}(G,\
  \Upsilon)^*\Big).
 \end{equation}
Here $v_2=\{v_{2}^{(k,\zeta)}\}$ is defined
by~(\ref{eqPowerSol2*}), $d_{2}=\{d_{2}^{(k,\zeta)}\}$ is a vector
of some constants; $v_{21}=\{v_{21}^{(k,\zeta)}\}$, where
$v_{21}^{(k,\zeta)}$ is a linear combination of the functions
$r^{i\bar\lambda_1+2m-2}(i\ln r)^q\Psi_{21q}^{(k,\zeta)}$, $0\le
q\le k+\varkappa(\bar\lambda_1)$, $\Psi_{21q}^{(k,\zeta)}\in
W_{2,2\pi}^l(0,\ 2\pi)^*$; $d_{1}$ is the vector of constants
appearing in~(\ref{eqAsymp*1}).
\end{theorem}
\begin{proof}
Let $v_{21}=\{v_{21}^{(k,\zeta)}\}$ be a particular solution
(which is defined by Lemma~5.2~\cite{GurPlane}) for the problem
\begin{equation}\label{eqL20Spec}
 {\cal L}_2^*v_{21}=-\Psi_{21},
\end{equation}
where $\Psi_{21}=\{\Psi_{21}^{(k,\zeta)}\}$ is a ``special''
distribution appearing in~(\ref{eqL*20}).
 We remind that each element
$\Psi_{21}^{(k,\zeta)}$ is a linear combination of the
distributions $r^{i\bar\lambda_{1}-2}(i\ln
r)^q\Psi_{21q}^{(k,\zeta)}$, $0\le q\le k$. Therefore, by
Lemma~5.2~\cite{GurPlane}, the particular solution $v_{21}$ has
the form described in the formulation of the theorem. Moreover,
each component $v_{21}^{(k,\zeta)}$ of the vector $v_{21}$ is
uniquely defined if $\bar\lambda_1$ is not an eigenvalue of
$\tilde{\cal L}^*_2(\lambda)$ (i.e., if
$\bar\lambda_1\ne\bar\lambda_2$). Otherwise (i.e., if
$\bar\lambda_1=\bar\lambda_2$) it is defined accurate to an
arbitrary linear combination of power
solutions~(\ref{eqPowerSol2*}) corresponding to the
eigenvalue~$\bar\lambda_1=\bar\lambda_2$. From now on we shall
suppose a particular solution $v_{21}=\{v_{21}^{(k,\zeta)}\}$
being fixed.

Combining~(\ref{eqL*20}) with~(\ref{eqL20Spec}) and using
Leibnitz's formula, one easily checks that ${\cal L}_2^*
(\eta_2v-d_1\eta_2v_{21})\in H_{a}^{l+2m}({\mathbb R}^2)^*$. Now
the asymptotics~(\ref{eqAsymp*2}) is resulted from
Theorem~5.3~\cite{GurPlane}, which establishes the asymptotics of
solutions for adjoint problems in ${\mathbb R}^2$.
\end{proof}

Theorem~\ref{thAsymp*2} shows that the asymptotic behavior of
solutions for adjoint nonlocal problem~(\ref{eqL*G}) near the
point $g_2$ depends on the data of the problem both near the point
$g_2$ itself and near the point $g_1$, which is connected with
$g_2$: $g_1=\Omega_1^{-1}(g_2)$.

\medskip

{\bf IV.} Let us write the asymptotics formula for the solution
$\{v,\ w_{\sigma}\}\in H_{a_1}^l(G,\ \Upsilon)^*$ for adjoint
nonlocal problem~(\ref{eqL*G}) in the whole domain $G$. Put
(cf.~(\ref{eqU12}))
\begin{equation}\label{eqV21}
 \{V_{2},\ W_{2,\sigma}\}=\eta_2 \{v_{2},\ 0\};\quad
 \{V_{1},\ W_{1,\sigma}\}=\eta_1 \{v_{1},\ w_{1,\sigma}\}
 + \eta_2 \{v_{21},\ 0\}.
\end{equation}

Now Theorems~\ref{thAsymp*1} and~\ref{thAsymp*2} yield the
following asymptotics of $\{v,\ w_{\sigma}\}\in H_{a_1}^l(G,\
\Upsilon)^*$ (cf. formula~(\ref{eqAsympG})):
\begin{equation}\label{eqAsymp*G}
\{v,\ w_{\sigma}\}\equiv \Big(d_1 \{V_{1},\ W_{1,\sigma}\} +
d_2\{V_{2},\ W_{2,\sigma}\}\Big) \Big({\rm mod\,} H_{a}^{l+2m}(G,\
\Upsilon)^*\Big).
\end{equation}

\section{Calculation of the coefficients
in the asymptotics formulas}
\phantom{a}~\newline\nopagebreak\label{sectCoef}

{\bf I.} In this section we will calculate the coefficients
$c_{\nu}^{(k,\zeta)}$ appearing in asymptotics~(\ref{eqAsympG}).

To begin with, let us remark that the coefficients can be
calculated in the following way. At first one should find
$c_{1}^{(k,\zeta)}$. Since in the neighborhood ${\cal V}(g_2)$ of
the point $g_2$ the function $u$ has asymptotics~(\ref{eqAsymp2}),
by Theorem~5.2~\cite{GurPlane} we have
\begin{equation}\label{eqc2}
c_{2}^{(k,\zeta)}=<{\cal L}_2 \eta_2 u,\ i
v_{2}^{(\varkappa_{\zeta,2}-k-1,\zeta)}>,
\end{equation}
where $v_{2}^{(k,\zeta)}$ is defined in~(\ref{eqPowerSol2*}).
Further, by Theorem~\ref{thAsymp1}, the function $u'=u-c_2u_{12}$
(where $c_2$ is calculated in~(\ref{eqc2}), $u_{12}$ is defined in
the proof of Theorem~\ref{thAsymp1}) has the following asymptotics
in the neighborhood ${\cal V}(g_1)$ of the point $g_1$:
\begin{equation}\label{equ'}
 u'(y)=c_1u_1(y)+\hat u(y)\quad (y\in {\cal V}(g_1)\cap G).
\end{equation}
Here $u_1$ is defined by~(\ref{eqPowerSolLambda_n1}); $c_1$ is to
be found; $\hat u\in H_{a_1}^{l+2m}({\cal V}(g_1)\cap G)$. From
asymptotics~(\ref{equ'}) and Theorem~4.1~\cite{GurPlane}, it
follows that
\begin{equation}\label{eqc1}
c_{1}^{(k,\zeta)}=<{\cal L}_1 \eta_1 u',\ i
\{v_{1}^{(\varkappa_{\zeta,1}-k-1,\zeta)},\
w_{1,\sigma}^{(\varkappa_{\zeta,1}-k-1,\zeta)}\}>,
\end{equation}
where $\{v_{1}^{k,\zeta)},\ w_{1,\sigma}^{(k,\zeta)}\}$ is defined
in~(\ref{eqPowerSol1*}).

Formulas~(\ref{eqc2}) and~(\ref{eqc1}) show that the value of
$c_1$ (as well as the general form of the asymptotics near $g_1$)
depends not only on the data of the problem near the point $g_1$
but also from the data near $g_2=\Omega_1(g_1)$.

\smallskip

We remark that similarly to~(\ref{eqc2}) and~(\ref{eqc1}) one can
calculate the coefficients $c_\nu$ with the help of the Green
formula and so--called formally adjoint problems generated by the
Green formula\footnote{In this case, additionally to
Conditions~\ref{condEllipPinG} and~\ref{condComplBinG}, one must
demand the system $\{B_{\sigma\mu}(D_y)\}_{\mu=1}^m$ to be normal
on $\gamma_\sigma$ ($\sigma=1,\ 2$), where $B_{\sigma\mu}(D_y)$ is
the principle homogeneous part of $B_{\sigma\mu}(g_1,\ D_y)$.}.
The corresponding technique is developed in~\cite{GurGiess,
GurPlane}. We will not recall the Green formula here, but only
mention that the corresponding formulas for $c_\nu$ are
immediately obtained if we use Theorems~5.4~\cite{GurPlane}
and~4.3~\cite{GurPlane} instead of Theorems~5.2~\cite{GurPlane}
and~4.1~\cite{GurPlane} respectively. Formally adjoint problems
have the advantage that they are considered in ``original''
spaces, but not in adjoint ones (spaces of distributions).
Therefore corresponding eigenvectors and associated vectors can be
found explicitly in a number of cases.

But, anyway, both adjoint problem-- and formally adjoint
problem--based formulas for $c_\nu$ involve the solution $u$
itself. Further we are to get formulas allowing to calculate the
coefficients $c_\nu$ only in terms of a right--hand side $\{f,\
f_{\sigma}\}$ of problem~(\ref{eqPinG}), (\ref{eqBinG}).

\medskip

{\bf II.} We are supposed to calculate $c_\nu$ with the help of
some special distributions from the kernel of the operator ${\bf
L}^*:H_{a_1}^l(G,\ \Upsilon)^*\to H_{a_1}^{l+2m}(G)^*$. To begin
with, assume that $\{v,\ w_\sigma\}\in H_{a_1}^l(G,\ \Upsilon)$ is
an arbitrary distribution from the kernel of ${\bf L}^*$.

 Let us calculate the value of the expression  $<{\bf L}u,\
i\{v,\ w_{\sigma}\}>$.

We suppose that the following consistent condition is fulfilled.
If the vector $c_\nu$ contains $c_{\nu}^{(k,\zeta)}$ in its $t$th
position, then the vector $d_\nu$ has
$d_{\nu}^{(\varkappa_{\zeta,\nu}-k-1,\zeta)}$ in its $t$th
position. The same is true for all the other vectors related to
the adjoint problem ($\{v_1,\ w_{1,\sigma}\}$, $v_2$, etc.).

Besides, we keep assuming that the Jordan chains corresponding to
the eigenvalues $\lambda_\nu$ and $\bar\lambda_\nu$ of the
operators $\tilde{\cal L}_\nu(\lambda)$ and $\tilde{\cal
L}^*_\nu(\lambda)$ respectively satisfy conditions of
biorthogonality and normalization~(\ref{eqNorm1}) (for $\nu=1$)
and~(\ref{eqNorm2}) (for $\nu=2$).

\smallskip

By virtue of~(\ref{eqAsympG}), we have ${\bf L}u={\bf L}
c_1U_1+{\bf L}c_2U_2+{\bf L}\hat u$, where $\hat u\in
H_{a_1}^{l+2m}(G)$. Since $\{v,\ w_{\sigma}\}$ belongs to the
kernel of ${\bf L}^*:H_{a_1}^l(G,\ \Upsilon)^*\to
H_{a_1}^{l+2m}(G)^*$, we get $<{\bf L}\hat u,\ i\{v,\
w_{\sigma}\}>=0$. Therefore, by virtue of~(\ref{eqSmoothLU}), we
can write
\begin{equation}\label{eqLuv}
<{\bf L}u,\ i\{v,\ w_{\sigma}\}>=<{\bf L}c_1 U_1,\ i\{v,\
w_{\sigma}\}>+<{\bf L}c_2 U_2,\ i\{v,\ w_{\sigma}\}>.
\end{equation}

Let $\eta_\nu(\omega,\ r)$ be the function $\eta_\nu$ written in
polar coordinates with the pole at $g_\nu$ ($\nu=1,\ 2$). For
$\varepsilon>0$ we introduce the functions
$\eta_{\nu,\varepsilon}(\omega,\ r)=\eta_\nu(\omega,\
r/\varepsilon)$.

At first let us consider the 1st term in the right--hand side
of~(\ref{eqLuv}). Since the difference
$\eta_{1}-\eta_{1,\varepsilon}$ vanishes near $g_1$, we have
$(\eta_{1}-\eta_{1,\varepsilon})c_1u_1\in H_{a_1}^{l+2m}(G)$. It
follows from this and from~(\ref{eqU12}) that
\begin{equation}\label{eqLuv1}
<{\bf L}c_1 U_1,\ i\{v,\ w_\sigma\}>=<{\bf L}c_1
\eta_{1,\varepsilon}u_1,\ i\{v,\ w_{\sigma}\}>.
\end{equation}

Put for short ${\cal U}_{1,\varepsilon}=c_1
\eta_{1,\varepsilon}u_1$. Since ${\rm supp\,}{\cal
U}_{1,\varepsilon}\subset{\cal V}(g_1)\cap G={\cal V}(0)\cap K$,
we have
$$
{\bf L}{\cal U}_{1,\varepsilon}= {\cal L}{\cal
U}_{1,\varepsilon}+\{{\bf P}(y,\ D_y)-{\cal P}(D_y),\ {\bf
B}_\sigma(y,\ D_y)-{\cal B}_\sigma(D_y)\}{\cal U}_{1,\varepsilon}.
$$
Here ${\cal P}(D_y)$ is the principal homogeneous part of ${\bf
P}(g_1,\ D_y)$; ${\cal B}_\sigma(D_y)$ is defined
by~(\ref{eqcalB1}). This and Theorem~\ref{thAsymp*1} imply that
the right--hand side of~(\ref{eqLuv1}) has the form
\begin{multline}\label{eqLuv1'}
<{\cal L}_1{\cal U}_{1,\varepsilon},\ i d_1\{v_1,\ w_{1,\sigma}\}>
+<{\cal L}_1{\cal U}_{1,\varepsilon},\ i\{F,\ G_\sigma\}>+\\
 <\{{\bf P}(y,\ D_y)-{\cal P}(D_y),\  {\bf B}_{\sigma}(y,\ D_{y})
- {\cal B}_{\sigma}(D_y)\} {\cal U}_{1,\varepsilon},\\
    i (d_1\{v_1,\ w_{1,\sigma}\}+\{F,\ G_{\sigma}\})>,
\end{multline}
where $\{v_1,\ w_{1,\sigma}\}$ is defined by~(\ref{eqPowerSol1*}),
$\{F,\ G_{\sigma}\}\in H_a^{l}(K,\ \gamma)^*$.

By Theorem~4.1~\cite{GurPlane}, the 1st term in~(\ref{eqLuv1'}) is
equal to
 $(c_1,\
d_1)$\footnote{Here $(c_1,\ d_1)$ (and further $(c_2,\ d_2)$, etc)
stands for the inner product of the corresponding complex vectors:
$(c_1,\ d_1)=\sum\limits_{k,\zeta}c_{1}^{(k,\zeta)} \overline
{d_{1}^{(\varkappa_{\zeta,1}-k-1,\zeta)}}$.}. The 2nd term
in~(\ref{eqLuv1'}) is majorized by
$$
c\|{\cal U}_{1,\varepsilon}\|_{H_a^{l+2m}(K)} \|\{F,\
G_{\sigma}\}\|_{H_a^{l}(K,\ \gamma)^*}=
 {\rm O}(1),
$$
where we use the Hardy--Littlewood symbol ``O'' with its usual
interpretation (${\rm O}(1)$ tends to $0$ as $\varepsilon\to0$).

By virtue of the boundedness of the imbedding operator of
$H_{a_1+1}^{l+2m}(K)$ into $H_{a_1}^{l+2m-1}(K)$,
Lemma~3.3$'$~\cite{KondrTMMO67}, and the inequality $a<a_1+1$, the
last term in~(\ref{eqLuv1'}) is majorized by
$$
c\|{\cal U}_{1,\varepsilon}\|_{H_{a_1+1}^{l+2m}(K)} \|d_1\{v_1,\
w_{1,\sigma}\}+\{F,\ G_{\sigma}\}\|_{H_{a_1}(K,\ \gamma)^*}\le c'
 \|{\cal U}_{1,\varepsilon}\|_{H_{a}^{l+2m}(K)} ={\rm O}(1).
$$
Thus, as $\varepsilon$ tends to 0, we get
\begin{equation}\label{eqLuvNu1}
 <{\bf L}c_1 U_1,\ i\{v,\ w_{\sigma}\}>=(c_1,\ d_1).
\end{equation}

\smallskip

Now let us consider the 2nd term in the right--hand side
of~(\ref{eqLuv}). Since the functions
$(\eta_{2}-\eta_{2,\varepsilon})c_2u_2$ and
$(\eta_{1}-\eta_{1,\varepsilon})c_2u_{12}$ belong to the space
$H_{a_1}^{l+2m}(G)$, we obtain from~(\ref{eqU12}) that
\begin{equation}\label{eqLuv2}
<{\bf L}c_2 U_2,\ i\{v,\ w_\sigma\}>=<{\bf L}c_2(
\eta_{2,\varepsilon}u_2+\eta_{1,\varepsilon}u_{12}),\ i\{v,\
w_{\sigma}\}>.
\end{equation}

Put for short ${\cal U}_{2,\varepsilon}=c_2
\eta_{2,\varepsilon}u_2$, ${\cal U}_{12,\varepsilon}=c_2
\eta_{1,\varepsilon}u_{12}$. Using Theorems~\ref{thAsymp*1}
and~\ref{thAsymp*2}, write the right--hand side of~(\ref{eqLuv2})
in the form
\begin{multline}\label{eqLuv2'}
<{\cal L}_2{\cal U}_{2,\varepsilon},\ i d_2v_2> +<{\cal L}_2{\cal
U}_{2,\varepsilon},\ i d_1v_{21}> +<{\cal P}(D_y){\cal
U}_{12,\varepsilon},\ id_1v_1>+\\
<{\cal B}_{1}(D_y) {\cal U}_{12,\varepsilon}+\big(T_1(D_y){\cal
U}_{2,\varepsilon}\big)\big(\Omega_1(y)\big)|_{\gamma_1},\ i d_1
w_{1,1}>+\\
<{\cal B}_{2}(D_y) {\cal U}_{12,\varepsilon}\ i d_1
w_{1,2}>+{\rm O}(1).
\end{multline}
Here ${\cal P}(D_y),$ $T_1(D_y)$ are the principal homogeneous
parts of ${\bf P}(g_1,\ D_y)$, ${\bf T}(g_2,\ D_y)$ respectively;
${\cal B}_\sigma(D_y)$, $\sigma=1,\ 2$, are defined
by~(\ref{eqcalB1}).

By Theorem~5.2~\cite{GurPlane}, the 1st term in~(\ref{eqLuv2'}) is
equal to $(c_2,\ d_2)$.

Since ${\cal L}_2c_2u_2=0$ (see~(\ref{eqL2u2=0})), the 2nd term
in~(\ref{eqLuv2'}) is equal to
\begin{equation}\label{eqLuvNu22}
<[{\cal L}_2,\ \eta_{2,\varepsilon}] c_2u_2,\
     i d_1 v_{21}>,
\end{equation}
where $[\cdot,\ \cdot]$ is the commutator. Using the condition
$0<a-a_1<1$, one can easily check that~(\ref{eqLuvNu22}) is equal
to
\begin{equation}\label{eqLuvNu23}
(\hat A_{12}(\varepsilon) c_2,\ d_1),
\end{equation}
where $\hat A_{12}(\varepsilon)$ is a matrix of the corresponding
order, the elements of which are linear combinations of the
functions $\varepsilon^{\lambda_2-\lambda_1}(i\ln\varepsilon)^q$.

Further, let us recall that the function $u_{12}$ is a solution
for problem~(\ref{eqP10Spec}), (\ref{eqB10Spec}). Hence the sum of
the 3rd, 4th, and 5th terms in~(\ref{eqLuv2'}) is equal to
\begin{multline}\label{eqLuvNu22'}
<[{\cal P}(D_y),\ \eta_{1,\varepsilon}] c_2u_{12},\ i d_1
v_{1}>+<[{\cal B}_{1}(D_y),\
\eta_{1,\varepsilon}]c_2u_{12}+\\
\big([T_1(D_y),\
\eta_{1,\varepsilon}]c_2u_{2}\big)\big(\Omega_1(y)\big)|_{\gamma_1},\
i d_1 w_{1,1}>+\\
<[{\cal B}_{2}(D_y),\
\eta_{1,\varepsilon}]c_2u_{12},\ i d_1 w_{1,2}>
\end{multline}
and therefore is of the form~(\ref{eqLuvNu23}). Thus we see that
\begin{equation}\label{eqLuvNu2}
<{\bf L}c_2 U_2,\ i\{v,\ w_{\sigma}\}>=(c_2,\ d_2)+
 (A_{12}(\varepsilon) c_2,\ d_1)+{\rm O}(1),
\end{equation}
where $A_{12}(\varepsilon)$ is a matrix of the corresponding
order, the elements of which are linear combinations of the
functions $\varepsilon^{\lambda_2-\lambda_1}(i\ln\varepsilon)^q$.

From equations~(\ref{eqLuv}), (\ref{eqLuvNu1}),
and~(\ref{eqLuvNu2}), it follows that
\begin{equation}\label{eqLuvCoef}
<{\bf L}u,\ i\{v,\ w_{\sigma}\}>=(c_2,\ d_2)+
 (c_1+A_{12}(\varepsilon) c_2,\ d_1)+{\rm O}(1),
\end{equation}

\medskip

{\bf III.} Keeping denotation of section~\ref{sectIndex}, we will
denote by ${\cal U}_1,\ \dots,\ {\cal U}_\varkappa$ the ordered
set of functions $U_{\nu}^{(k,\zeta)}$, which are the elements of
the vectors $U_\nu$, $\nu=1,\ 2$, defined by~(\ref{eqU12}).

Denote by ${\cal V}_{1},\ \dots, {\cal V}_{\varkappa}$ the set of
functions $\{V_{\nu}^{(k,\zeta)},\ W_{\nu,\sigma}^{(k,\zeta)}\}$,
which are the elements of the vectors $\{V_{\nu},\
W_{\nu,\sigma}\}$, $\nu=1,\ 2$, defined by~(\ref{eqV21}).

Suppose the sets $\{{\cal U}_1,\ \dots, {\cal U}_\varkappa\}$ and
$\{{\cal V}_1,\ \dots, {\cal V}_\varkappa\}$ are ordered
consistently, i.e., the equality ${\cal U}_t=U_{\nu}^{(k,\zeta)}$
is fulfilled simultaneously with the equality ${\cal
V}_t=\{V_{\nu}^{(\varkappa_{\zeta,\nu}-k-1,\zeta)},\
W_{\nu,\sigma}^{(\varkappa_{\zeta,\nu}-k-1,\zeta)}\}$.

In this work we restrict ourselves to the case when $d=0$ in
Lemma~\ref{lKerL*}. This mean that any solution to homogeneous
problem~(\ref{eqPinG}), (\ref{eqBinG}) from the space
$H_a^{l+2m}(G)$ necessarily belongs to the space
$H_{a_1}^{l+2m}(G)$. In that case we will show that for any
right--hand side $\{f,\ f_{\sigma}\}\in H_{a_1}^{l}(G,\ \Upsilon)$
the coefficients in the asymptotics formula for solutions are
uniquely defined. If $d>0$, then, similarly to the case of
``local'' problems (see Theorem~3.6~\cite[Chapter~4]{NP}), there
is some freedom in choosing the coefficients of the asymptotics.
Moreover, the procedure for calculation of the coefficients
becomes more technically complicated (while the idea remains
similar to the one we shall describe below) and will not be
considered here.

So, suppose $d=0$. Then, by virtue of Lemma~\ref{lKerL*}, there
exist solutions ${\cal Y}_{1},\ \dots,\ {\cal Y}_\varkappa\in
H_{a_1}^l(G,\ \Upsilon)^*$ for the equation ${\bf L}^*{\cal Y}=0$,
linearly independent modulo $H_{a}^l(G,\ \Upsilon)^*$.
By~(\ref{eqAsymp*G}) we have ${\cal
Y}_t\equiv\sum\limits_{k=1}^{\varkappa}d_{tk}{\cal V}_k\ \Big({\rm
mod\,} H_{a}^l(G,\ \Upsilon)^*\Big)$, $t=1,\ \dots,\ \varkappa$.
Since ${\cal Y}_{1},\ \dots,\ {\cal Y}_\varkappa$ are linearly
independent modulo $H_{a}^l(G,\ \Upsilon)^*$, the matrix
$\|d_{tk}\|$ is nonsingular. Hence, without loss in generality, we
can assume that
\begin{equation}\label{eq_ExCanonicBasis*}
{\cal Y}_t\equiv {\cal V}_t\ \Bigl({\rm mod\,} H_{a}^l(G,\
\Upsilon)^*\Bigr),\quad t=1,\ \dots,\ \varkappa.
\end{equation}

Now let us prove that the elements of the matrix
$A_{12}(\varepsilon)$ appearing in~(\ref{eqLuvCoef}) have finite
limits as $\varepsilon\to0$. This limit will be denoted by
$A_{12}$:
$$
 A_{12}=\lim\limits_{\varepsilon\to0}A_{12}(\varepsilon).
$$
Let $l_\nu$ be the length of the vector $c_\nu$ (or $d_\nu$, which
is the same), $\nu=1,\ 2$. Clearly, $l_2+l_1=\varkappa$. Suppose
for definiteness that the first $l_2$ elements in the ordered set
$\{{\cal U}_1,\ \dots, {\cal U}_\varkappa\}$ ($\{{\cal V}_1,\
\dots, {\cal V}_\varkappa\}$) are components of the vector $U_2$
($\{V_2,\ W_{2,\sigma}\}$) and the last $l_1$ ones are components
of the vector $U_1$ ($\{V_1,\ W_{1,\sigma}\}$):
\begin{gather}
\{{\cal U}_1,\ \dots, {\cal
U}_\varkappa\}=\{\underbrace{U_2}_{l_2},\
\underbrace{U_1}_{l_1}\}\notag\\
 \Big(\{{\cal V}_1,\ \dots, {\cal
V}_\varkappa\}=\big\{\underbrace{\{V_2,\ W_{2,\sigma}\}}_{l_2},\
\underbrace{\{V_1,\ W_{1,\sigma}\}}_{l_1}\big\}\Big).\notag
\end{gather}

Now fix an arbitrary $t$ from the set $\{1,\ \dots,\ l_2\}$ and an
arbitrary $k$ from the set $\{l_2+1,\ \dots,\ \varkappa\}$.
Substituting in~(\ref{eqLuvCoef}) $u={\cal U}_{t}$ (which is a
component of the vector $U_2$) and $\{v,\ w_{\sigma}\}={\cal
Y}_{k}$ (which is, by~(\ref{eq_ExCanonicBasis*}), a component of
the vector $\{V_1,\ W_{1,\sigma}\}$ modulo $H_{a}^l(G,\
\Upsilon)^*$), we get $c_1=0$, $d_2=0$, and therefore,
\begin{equation}\label{eqa_tk}
 <{\bf L}U_{t},\ i{\cal Y}_{k}\}>=a_{tk}(\varepsilon)+{\rm O}(1).
\end{equation}
Here $a_{tk}(\varepsilon)$ is the corresponding element of the
matrix $A_{12}(\varepsilon)$. The left--hand side
of~(\ref{eqa_tk}) does not depend on $\varepsilon$. Therefore
$a_{tk}(\varepsilon)$ has a finite limit as $\varepsilon\to0$.

Thus, passing in~(\ref{eqLuvCoef}) to the limit as
$\varepsilon\to0$, we get
\begin{equation}\label{eqLuvCoefConst}
<{\bf L}u,\ i\{v,\ w_{\sigma}\}>=(c_2,\ d_2)+
 (c_1+A_{12} c_2,\ d_1).
\end{equation}

\begin{theorem}\label{thCoef}
Let $u\in H_a^{l+2m}(G)$ be a solution for problem~(\ref{eqPinG}),
(\ref{eqBinG}) with a right--hand side $\{f,\ f_{\sigma}\}\in
H_{a_1}^l(G,\ \Upsilon)$. Then $u$ has the asymptotics
 \begin{equation}\label{eqCoef}
  u\equiv \Big(\sum\limits_{t=1}^\varkappa c_t {\cal U}_t\Big)\ \Big({\rm mod\,} H_{a_1}^{l+2m}(G)\Big).
 \end{equation}
The constants $c_t$ ($t=1,\ \dots,\ \varkappa$) can be calculated
by the formulas
 \begin{equation}\label{eqCoef2}
  c_t=<\{f,\ f_{\sigma}\},\ i{\cal Y}_{t}>
 \end{equation}
if $t\le l_2$ (i.e., $c_t$ coincides with a component of the
vector $\{c_{2}^{(k,\zeta)}\}$);
 \begin{equation}\label{eqCoef1}
  c_t=<\{f,\ f_{\sigma}\},\ i{\cal Y}_{t}-
  i \Big[A_{12}({\cal Y}_{1},\ \dots,\ {\cal Y}_{l_2})^T\Big]_{t-l_2}>
 \end{equation}
if $l_2<t\le \varkappa$ (i.e., $c_t$ coincides with a component of
the vector $\{c_{1}^{(k,\zeta)}\}$). Here $[\cdot]_j$ stands for
the $j$th component of a vector.
\end{theorem}
\begin{proof}
Substituting $\{v,\ w_{\sigma}\}={\cal Y}_{1},\ \dots,\ \{v,\
w_{\sigma}\}={\cal Y}_{\varkappa}$ subsequently
in~(\ref{eqLuvCoefConst}), we obtain formulas~(\ref{eqCoef2})
and~(\ref{eqCoef1}).
\end{proof}

Theorem~\ref{thCoef} shows that the values of the coefficients
$c_{\nu}^{(k,\zeta)}$ are the functionals over the right--hand
sides $\{f,\ f_{\sigma}\}$ of problem~(\ref{eqPinG}),
(\ref{eqBinG}). These functionals depend on the data of the
problem in the whole domain $G$, but not only in the neighborhoods
${\cal V}(g_1)$ and ${\cal V}(g_2)$.

\begin{remark}
We remind that the elements of the matrix $A_{12}(\varepsilon)$
are linear combinations of the functions
$\varepsilon^{\lambda_2-\lambda_1}(i\ln\varepsilon)^q$. Hence, if
$\lambda_1\ne\lambda_2$, then $A_{12}=0$.
\end{remark}

\section{Example}
\phantom{a}~\newline \nopagebreak\label{sectExample}

{\bf I.} In this section we consider an example illustrating the
results of sections~\ref{sectStatement}--\ref{sectCoef}.

Keeping denotation and assumptions of sections~\ref{sectStatement}
and~\ref{sectAsymp}, we consider the following nonlocal problem
\begin{gather}
{\bf P}(y,\
D_y)\equiv\sum_{|\alpha|\le2}p_\alpha(y)\frac{\partial^{|\alpha|}u}{\partial
y_1^{\alpha_1}\partial y_2^{\alpha_2}}=f(y) \quad (y\in
G\setminus{{\cal
 K}}),\label{eq_ExPinG}\\
{\bf B}_{\sigma }u\equiv u(y)|_{\Upsilon_\sigma}+ e_\sigma
u\big(\Omega_{\sigma }(y)\big)|_{\Upsilon_\sigma}=f_{\sigma
}(y)\quad
    (y\in \Upsilon_\sigma;\ \sigma=1,\ 2)\label{eq_ExBinG}.
\end{gather}
Here ${\bf P}(y,\ D_y)$ is a 2nd order differential operator,
properly elliptic in $\bar G$, with infinitely smooth coefficients
$p_\alpha(y)$; $e_\sigma\in\mathbb C$. For clearness we assume
\begin{equation}\label{eqPDelta}
 \sum_{|\alpha|=2}p_\alpha(g_\nu)\frac{\partial^{|\alpha|}u}{\partial
y_1^{\alpha_1}\partial y_2^{\alpha_2}}=\Delta u,\quad \nu=1,\ 2.
\end{equation}

\smallskip

Let us obtain the asymptotics of a solution $u\in H_a^{2}(G)$ for
problem~(\ref{eq_ExPinG}), (\ref{eq_ExBinG}) with a right--hand
side $\{f,\ f_{\sigma }\}\in H_{a_1}^0(G,\
\Upsilon)\stackrel{def}{=}H_{a_1}^0(G)\times\prod\limits_{\sigma=1,2}
H_{a_1}^{3/2}(\Upsilon_\sigma)$, assuming $0<a-a_1<1$.

At first, according to section~\ref{sectAsymp}, we consider the
asymptotics of the solution $u$ in the neighborhood ${\cal
V}(g_2)$ of the point $g_2$. For this purpose one must write the
model equation in ${\mathbb R}^2\setminus\{g_2\}$. Taking into
account~(\ref{eqPDelta}), we obtain
\begin{equation}\label{eq_ExP3}
 \Delta u=\hat f(y) \quad (y\in{\mathbb R}^2\setminus\{g_2\}),
\end{equation}
where $\hat f\in H_{a_1}^0\big({\cal V}(g_2)\big)$.

Write equation~(\ref{eq_ExP3}) in polar coordinates with the pole
at $g_2$:
$$
 r\frac{\partial}{\partial r}\Big(r\frac{\partial u}{\partial
 r}\Big)+\frac{\partial^2u}{\partial \omega^2}=r^2f(\omega,\
 r)\quad (0<\omega<2\pi,\ r>0).
$$
Applying formally the Mellin transformation, we get
$$
\frac{d^2\tilde u}{d \omega^2}-\lambda^2\tilde u=\tilde
F(\lambda,\ \omega)\quad (0<\omega<2\pi),
$$
where $\tilde u$ and $\tilde F$ are the Mellin transforms of $u$
and $r^2f$ with respect to r.

Introduce the corresponding operator--valued function
$$
  \tilde{\cal
L}_2(\lambda)=\frac{\displaystyle d^2}{\displaystyle
d\omega^2}-\lambda^2 :
 W_{2\pi}^{2}(0,\ 2\pi)\to L_{2}(0,\ 2\pi).
$$

Let us suppose, additionally to Condition~\ref{condSpectr}, that
there is the only eigenvector $\varphi_2(\omega)$ corresponding to
the eigenvalue $\lambda_2$ ($a_1-1<{\rm Im\,}\lambda_2<a-1$) of
$\tilde{\cal L}_2(\lambda)$ and there are no associated vectors.

Then, by Theorem~\ref{thAsymp2}, we have
 \begin{equation}\label{eq_ExAsymp2}
  u(y)=c_2 u_2(y)+\hat u(y)\quad \big(y\in{\cal V}(g_2)\big).
 \end{equation}
Here $c_2$ is a scalar constant,
$u_2=r^{i\lambda_{2}}\varphi_{2}(\omega)$ is a power solution for
homogeneous equation~(\ref{eq_ExP3}); $(\omega,\ r)$ are polar
coordinates with the pole at $g_2$ and the polar axis being, for
definiteness, tangent to the curve $\Omega_1(\Upsilon_1)$ at
$g_2$; $\hat u\in H_{a_1}^{2}\big({\cal V}(g_2)\big)$.

\smallskip

Now we consider the asymptotics of the solution $u$ for
problem~(\ref{eq_ExPinG}), (\ref{eq_ExBinG}) in the neighborhood
${\cal V}(g_1)$ of the point $g_1$. Let $\Omega_1(y)$ ($y\in{\cal
V}(g_1)$) be a rotation with respect to $g_1$ (with no expansion
for simplicity) and the shift by the vector
$\overrightarrow{g_1g_2}$. Let $\Omega_2(y)$ ($y\in{\cal V}(g_1)$)
coincide with the operator ${\cal G}_2$ of a rotation by an angle
$\omega_2$ ($b_1<b_2+\omega_2<b_2$) and an expansion with a
coefficient $\beta_2>0$.

According to section~\ref{sectAsymp} and
assumption~(\ref{eqPDelta}), the asymptotics of $u$ in ${\cal
V}(g_1)$ coincides with the asymptotics of a solution for the
problem
\begin{gather}
 \Delta u=\hat f(y) \quad (y\in{\cal V}(0)\cap K),\label{eq_ExP1}\\
\begin{aligned}\label{eq_ExB1}
u|_{ {\cal V}(0)\cap\gamma_1}=\hat f_1-c_2 f_{12}\quad
    (y\in {\cal V}(0)\cap\gamma_1),\\
u|_{ {\cal V}(0)\cap\gamma_2}+e_2 u({\cal G}_2y)|_{ {\cal
V}(0)\cap\gamma_2}=f_2\quad     (y\in {\cal V}(0)\cap\gamma_2).
\end{aligned}
\end{gather}
Here $\hat f\in H_{a_1}^0({\cal V}(0)\cap K)$,
\begin{align*}
 \hat f_1&=f_1-e_1\hat u\big(\Omega_{1}(y)\big)|_{ {\cal
V}(0)\cap\gamma_1}\in H_{a_1}^{3/2}({\cal
V}(0)\cap\gamma_1),\\
f_{12}&=e_1u_2\big(\Omega_{1}(y)\big)|_{ {\cal
V}(0)\cap\gamma_1}=e_1 r^{i\lambda_2}\varphi_2(0)\footnotemark{}.
\end{align*}
\footnotetext{We calculate $\varphi_2(\omega)$ for $\omega=0$
because of the special choice of polar coordinates (see above).}

Similarly to the above we obtain the corresponding
operator--valued function $ \tilde{\cal
L}_1(\lambda):W_2^{2}(b_1,\ b_2)\to L_2(b_1,\ b_2) \times{\mathbb
C}^{2}$
 given by
\begin{equation}\label{1_eqtildeLLambda}
\tilde{\cal L}_1(\lambda)\varphi=\Big\{\frac{\displaystyle d^2
\varphi }{\displaystyle d\omega^2} -\lambda^2 \varphi,\
\varphi(\omega)|_{\omega=b_1},\
 \varphi(\omega)|_{\omega=b_2}+e_2 e^{i\lambda\ln\beta_2 }
\varphi(\omega+\omega_2)|_{\omega=b_2}\Big\}.
\end{equation}

Let us suppose, additionally to Condition~\ref{condSpectr}, that
there is the only eigenvector $\varphi_1(\omega)$ corresponding to
the eigenvalue $\lambda_1$ ($a_1-1<{\rm Im\,}\lambda_1<a-1$) of
$\tilde{\cal L}_1(\lambda)$ and there are no associated vectors.

Then, by Theorem~\ref{eqAsymp1}, we have
 \begin{equation}\label{eq_ExAsymp1}
  u(y)=c_1 u_1(y)+c_2u_{12}(y)+\hat u(y)\quad (y\in{\cal V}(g_1)).
 \end{equation}
Here $c_1$ is some scalar constant, $c_2$ is the constant
appearing in~(\ref{eq_ExAsymp2});
$u_{12}=r^{i\lambda_2}\varphi_{12}(\omega)$ ($\varphi_{12}\in
W_2^2(b_1,\ b_2)$) is a particular solution for the following
problem in the angle $K$ with the ``special'' right--hand side
(cf.~(\ref{eqP10Spec}), (\ref{eqB10Spec})):
\begin{gather}
 \Delta u=0 \quad (y\in K),\label{eq_ExP1Spec}\\
 u|_{\gamma_1}=-f_{12},\quad u|_{\gamma_2}+e_2
u({\cal G}_2y)|_{\gamma_2}=0;\label{eq_ExB1Spec}
\end{gather}
$u_1=r^{i\lambda_{1}}\varphi_{1}(\omega)$ is a solution for
homogeneous problem~(\ref{eq_ExP1Spec}), (\ref{eq_ExB1Spec});
$(\omega,\ r)$ are polar coordinates with the pole at the point
$g_1=0$; $\hat u\in H_{a_1}^{2}\big({\cal V}(g_1)\big)$.

\smallskip

To write the asymptotics in the whole domain $G$, we introduce the
functions
$
 U_{1}=\eta_1 u_{1},\
 U_{2}=\eta_2 u_{2}+\eta_1 u_{12}.
$ Then~(\ref{eq_ExAsymp2}) and~(\ref{eq_ExAsymp1}) imply:

{\it Let $u\in H_a^{2}(G)$ be a solution for
problem~(\ref{eq_ExPinG}), (\ref{eq_ExBinG}) with a right--hand
side $\{f,\ f_{\sigma }\}\in H_{a_1}^0(G,\ \Upsilon)$,
$0<a-a_1<1$. Then we have
 \begin{equation}\label{eq_ExAsympG}
  u\equiv \Big(c_{1}U_{1}+c_{2}U_{2}\Big)\
 \Big({\rm mod\,} H_{a_1}^{2}(G)\Big),
 \end{equation}
where $c_1,\ c_2$ are some scalar constants. }

\medskip

{\bf II.} From asymptotics formula~(\ref{eq_ExAsympG}) and
Theorem~\ref{thIndL} we can derive the connection between the
indices of the operators
\begin{align*}
  {\bf L}_a&=\{{\bf P}(y,\ D_y),\ {\bf B}_{\sigma }\}: H_a^{2}(G)\to  H_a^0(G,\
  \Upsilon),\\
  {\bf L}_{a_1}&=\{{\bf P}(y,\ D_y),\ {\bf B}_{\sigma }\}: H_{a_1}^{2}(G)\to  H_{a_1}^0(G,\ \Upsilon)
\end{align*}
corresponding to problem~(\ref{eq_ExPinG}), (\ref{eq_ExBinG}), but
acting in different weighted spaces. Since the sum of full
multiplicities of eigenvalues $\lambda_1$ and $\lambda_2$ is equal
to 2 in our case, the connection between the indices is as
follows:
$$
   {\rm ind\,}{\bf L}_a={\rm ind\,}{\bf L}_{a_1}+2.
$$

\medskip

{\bf III.} To calculate the coefficients $c_{\nu}$ in
formula~(\ref{eq_ExAsympG}), we will study the asymptotics of
solutions for the adjoint nonlocal problem.

Consider the operator ${\bf L}^*:H_{a_1}^0(G,\ \Upsilon)^*\to
H_{a_1}^{2}(G)^*$, adjoint to ${\bf L}=\{{\bf P}(y,\ D_y),\ {\bf
B}_{\sigma }\}: H_{a_1}^{2}(G)\to H_{a_1}^0(G,\ \Upsilon)$. The
operator ${\bf L}^*$ is given by
$$
<u,\ {\bf L}^*\{v,\ w_{\sigma }\}>=<{\bf P}(y,\ D_y) u,\
v>_G+\sum\limits_{\sigma=1,2} <{\bf B}_{\sigma }u,\ w_{\sigma
}>_{\Upsilon_\sigma},
$$
where $\{v,\ w_{\sigma }\}\in H_{a_1}^0(G,\ \Upsilon)^*$, $u\in
H_{a_1}^2(G)$.

Let us study the asymptotics of a solution $\{v,\ w_{\sigma }\}\in
H_{a_1}^0(G,\ \Upsilon)^*$ for the problem
\begin{equation}\label{eq_ExL*G}
 {\bf L}^*\{v,\ w_{\sigma }\}=\Psi,
\end{equation}
where $\Psi\in H_a^{2}(G)^*$.

\smallskip

In~\cite{GurPlane} it is shown that $\bar\lambda_\nu$ is an
eigenvalue of the operator $\tilde{\cal L}_\nu^*(\lambda)$,
adjoint to $\tilde{\cal L}_\nu(\bar\lambda)$. Denote by
$\{\psi_1,\ \chi_{1,\sigma}\}\in L_2(b_1,\ b_2)\times{\mathbb
C}^2$ ($\psi_2\in L_2(0,\ 2\pi)$) the eigenvector of $\tilde{\cal
L}_1^*(\lambda)$ ($\tilde{\cal L}_2^*(\lambda)$) corresponding to
the eigenvalue $\bar\lambda_1$ ($\bar\lambda_2$). Conditions of
biorthogonality and normalization~(\ref{eqNorm1})
and~(\ref{eqNorm2}) assume the form
\begin{equation}\label{eq_ExNormal}
 <-2\lambda_\nu\varphi_\nu,\ \psi_\nu>=1\footnotemark{}.
\end{equation}
\footnotetext{One can show that $\lambda_\nu\ne0$ whenever there
are no associated vectors corresponding to the eigenvalue
$\lambda_\nu$. Hence there always exist vectors $\{\psi_1,\
\chi_{1,\sigma}\}$ and $\psi_2$ satisfying~(\ref{eq_ExNormal}).}
Put $
 \{v_{1},\ w_{1,\sigma}\}=
  \{r^{i\bar\lambda_{1}}\psi_{1},\ r^{i\bar\lambda_{1}-1}\chi_{1,\sigma}\}\quad
(v_{2}=r^{i\bar\lambda_{2}}\psi_{2}),
$
where $(\omega,\ r)$ are polar coordinates with the pole at $g_1$
(with the pole at $g_2$ and with the polar axis being tangent to
the curve $\Omega_1(\Upsilon_1)$ at $g_2$).

Further, by Theorem~\ref{thAsymp*1}, we have
 \begin{equation}\label{eq_ExAsymp*1}
   \eta_1\{v,\ w_{\sigma }\}\equiv
  d_{1}\eta_1\{v_{1},\ w_{1,\sigma}\}\ \Big({\rm mod\,} H_a^{0}(G,\ \Upsilon)^*\Big),
 \end{equation}
where $d_{1}$ is some scalar constant. By Theorem~\ref{thAsymp*2},
we have
\begin{equation}\label{eq_ExAsymp*3}
   \eta_2 v\equiv
\Big(  d_{2}\eta_2 v_{2}+d_1\eta_2 v_{21}\Big)\
 \Big({\rm mod\,} H_a^{0}(G)^*\Big).
\end{equation}
Here $d_{2}$ is some scalar constant, $d_1$ is the constant
appearing in~(\ref{eq_ExAsymp*1});
$v_{21}=r^{i\bar\lambda_{1}}\psi_{21}$; $(\omega,\ r)$ are polar
coordinates with the pole at $g_2$ and the polar axis being
tangent to the curve $\Omega_1(\Upsilon_1)$ at $g_2$;
$\psi_{21}\in L_{2}(0,\ 2\pi)$. Moreover, the distribution
$v_{21}$ is a particular solution for the following adjoint
equation in ${\mathbb R}^2\setminus\{g_2\}$ with the ``special''
right--hand side (cf.~(\ref{eqL20Spec})):
$$
 \int\limits_{{\mathbb R}^2}\Delta u\cdot\bar
v\,dy=\int\limits_0^\infty u(0,\ r)\cdot(\overline{-\bar
e_1\chi_{1,1} r^{i\bar\lambda_{1}-1}})dr \quad \mbox{for all }
u\in C_0^\infty({\mathbb R}^2\setminus\{g_2\})\footnotemark{}.
$$
\footnotetext{We calculate $u(\omega,\ r)$ for $\omega=0$ because
of the special choice of polar coordinates (see above).}

Put
$
 \{V_{2},\ W_{2,\sigma}\}=\eta_2 \{v_{2},\ 0\},\
  \{V_{1},\ W_{1,\sigma}\}=\eta_1 \{v_{1},\ w_{1,\sigma}\}
 + \eta_2 \{v_{21},\ 0\}.
$

Then~(\ref{eq_ExAsymp*1}) and~(\ref{eq_ExAsymp*3}) imply:

{\it Let $\{v,\ w_{\sigma }\}\in H_{a_1}^0(G,\ \Upsilon)^*$ be a
solution for problem~(\ref{eq_ExL*G}) with a right--hand side
$\Psi\in H_a^{2}(G)^*$. Then we have
 \begin{equation}\label{eq_ExAsymp*G}
  \{v,\ w_{\sigma }\}\equiv
 \Big(d_1\{V_{1},\ W_{1,\sigma}\}+d_2\{V_{2},\ W_{2,\sigma}\}\Big)\
 \Big({\rm mod\,} H_{a}^{2}(G,\ \Upsilon)^*\Big),
 \end{equation}
where $d_1,\ d_2$ are some constants. }

\medskip

{\bf IV.} Now let us calculate the coefficients $c_{\nu}$
appearing in~(\ref{eq_ExAsympG}). Formulas~(\ref{eqc2})
and~(\ref{eqc1}) assume the form
\begin{gather}
 c_{2}=<\Delta (\eta_2 u),\ i v_{2}>_{{\mathbb R}^2},\notag\\
c_{1}=<\{\Delta u',\ u'|_{\gamma_1},\ u'|_{\gamma_2}+e_2 u'({\cal
G}_2y)|_{\gamma_2}\},\ i\{v_{1},\ w_{1,\sigma}\}>,\notag
\end{gather}
where $u'=\eta_1(u-c_{2} u_{12})$.

\smallskip

Now let us write the formulas allowing to calculate the
coefficients $c_\nu$ only in terms of a right--hand side $\{f,\
f_{\sigma}\}$ of problem~(\ref{eq_ExPinG}), (\ref{eq_ExBinG})
(i.e., independent of a solution $u$).

Following section~\ref{sectCoef}, we assume for simplicity that
any solution to homogeneous problem~(\ref{eq_ExPinG}),
(\ref{eq_ExBinG}) from the space $H_a^{2}(G)$ necessarily belongs
to the space $H_{a_1}^{2}(G)$. Then there exist solutions ${\cal
Y}_{1},\ {\cal Y}_2\in H_{a_1}^0(G,\ \Upsilon)^*$ for the equation
${\bf L}^*{\cal Y}=0$, linearly independent modulo $H_{a}^0(G,\
\Upsilon)^*$ such that
$$
{\cal Y}_\nu\equiv \{V_\nu,\ W_{\nu,\sigma}\}\ \Bigl({\rm mod\,}
H_{a}^0(G,\ \Upsilon)^*\Bigr),\quad \nu=1,\ 2.
$$

Let $\eta_{\nu,\varepsilon}$ be the functions defined in
section~\ref{sectCoef}.

Then from Theorem~\ref{thCoef} we obtain the following result.

{\it Let $u\in H_a^{2}(G)$ be a solution for
problem~(\ref{eq_ExPinG}), (\ref{eq_ExBinG}) with a right--hand
side $\{f,\ f_{\sigma}\}\in H_{a_1}^0(G,\ \Upsilon)$. Then the
function $u\in H_a^{2}(G)$ has asymptotics~(\ref{eq_ExAsympG}).
The constants $c_\nu$ ($\nu=1,\ 2$) are calculated by the formulas
\begin{gather}
  c_2=<\{f,\ f_{\sigma}\},\ i{\cal Y}_1>,\notag\\
  c_1=<\{f,\ f_{\sigma}\},\ i({\cal Y}_1-
  A_{12}{\cal Y}_2)>.\notag
\end{gather}
Here $A_{12}$ is a scalar constant given by
\begin{multline}\label{eq_ExCoefC}
A_{12}=\lim\limits_{\varepsilon\to0}<\Delta(\eta_{2,\varepsilon}u_2),\
iv_{21}>+<\{\Delta(\eta_{1,\varepsilon}u_{12}),
\eta_{1,\varepsilon}u_{12}|_{\gamma_1}+\eta_{1,\varepsilon}f_{12}|_{\gamma_1},\\
\eta_{1,\varepsilon}u_{12}|_{\gamma_2}+e_2
(\eta_{1,\varepsilon}u_{12})({\cal G}_2y)|_{\gamma_2}\},\ i\{v_1,\
w_{1,\sigma}\}>,
\end{multline}
where the limit does exist.}

\begin{remark}
The function $u_2$ ($u_{12}$) is a solution for homogeneous
equation~(\ref{eq_ExP3}) (a solution for
problem~(\ref{eq_ExP1Spec}), (\ref{eq_ExB1Spec}) with the special
right--hand side $\{0,\ -f_{12},\ 0\}$). Therefore, similarly to
section~\ref{sectCoef}, one can easily check that
$
A_{12}=\lim\limits_{\varepsilon\to0}{\rm const}\cdot
\varepsilon^{i(\lambda_2-\lambda_1)}.
$
From this and from the existence of the limit
in~(\ref{eq_ExCoefC}), it follows that $A_{12}=0$ whenever
$\lambda_1\ne\lambda_2$.
\end{remark}

\bigskip

The author is grateful to Professor Alexander Skubachevskii for
attention to this work and valuable advice.


\begin{thebibliography}{99}

\bibitem{Antonev}
A.B. Antonevich, The index and the normal solvability of a general
elliptic boundary value problem with a finite group of
translations on the boundary, {\it Differentsial'nye Uravneniya,}
{\bf 8} (1972), 309-317; English transl. in {\it Differential
Equations, } {\bf 8} (1974).

\bibitem{Beals}
R. Beals, Nonlocal elliptic boundary value problems, {\it Bull.
Amer. Math. Soc.,} {\bf 70} (1964), 693-696.

\bibitem{BitzSam}
A.V. Bitsadze and A.A. Samarskii, On some simple generalizations
of linear elliptic boundary value problems, {\it Dokl. Akad. Nauk
SSSR,} {\bf 185} (1969), 739-740; English transl. in {\it Soviet
Math. Dokl.,} {\bf 10} (1969).

\bibitem{Browder}
F. Browder, Non-local elliptic boundary value problems, {\it Amer.
J. Math.,} {\bf 86} (1964), 735-750.


\bibitem{Carleman}
T. Carleman, Sur la th\'eorie des equations integrales et ses
applications, {\it Verhandlungen des Internat. Math. Kongr.
Z\"urich,} {\bf 1} (1932), 132-151.

\bibitem{ZhEid}
S.D. Eidel'man and N.V. Zhitarashu, Nonlocal boundary value
problems for elliptic equations, {\it Mat. Issled.,} {\bf 6}
(1971), 63-73 (Russian).

\bibitem{Feller}
W. Feller, Diffusion processes in one dimension, {\it Trans. Amer.
Math. Soc.,} {\bf 77} (1954), 1-30.

\bibitem{GS}
I.C. Gohberg and E.I. Sigal, An operator generalization of the
logarithmic residue theorem and the theorem of Rouch\'e, {\it Mat.
Sb.,} {\bf 84} ({\bf 126}) (1971), 607-629; English transl. in
{\it Math. USSR Sb.,} {\bf 13} (1971).


\bibitem{GurDAN}
P.L. Gurevich, Nonlocal elliptic problems in dihedral angles and
the Green formula {\it Dokl. Akad. Nauk,} {\bf 379} (2001),
735-738; English transl. in {\it Russian Acad. Sci. Dokl. Math.,}
(2001).

\bibitem{GurGiess}
P.L. Gurevich, Nonlocal problems for elliptic equations in
dihedral angles and the Green formula, {\it Mitteilungen aus dem
Math. Seminar Giessen, Math. Inst. Univ. Giessen, Germany,} {\bf
247} (2001), 1-74.

\bibitem{GurMatZam}
P.L. Gurevich, Solvability of nonlocal elliptic problems in
dihedral angles, {\it Mat. Zametki,} {\bf 72} (2002), 178-197;
English transl. in {\it Math. Notes,} {\bf 72} (2002).

\bibitem{GurPlane}
P.L. Gurevich, Asymptotics of solutions for nonlocal elliptic
problems in plane angles, {\it Tr. semin. im. I.G.~Petrovskogo,}
{\bf 23} (2003); Engilsh transl. in {\it J. Math. Sci., New York}
(2004).


\bibitem{GM}
A.K. Gushchin and V.P. Mikhailov, On solvability of nonlocal
problems for elliptic equations of second order, {\it Mat. sb.,}
{\bf 185} (1994), 121-160; English transl. in {\it Math. Sb.,}
(1994).

\bibitem{Kishk}
K.Yu. Kishkis, The index of a Bitsadze--Samarskii problem for
harmonic functions, {\it Differentsial'nye Uravneniya,} {\bf 24}
(1988), 105-110; English transl. in {\it Differential Equations,}
{\bf 24} (1988), 83--87.


\bibitem{KondrTMMO67}
V.A. Kondrat'ev, Boundary value problems for elliptic equations in
domains with conical or angular points, {\it Trudy Moskov. Mat.
Obshch.,} {\bf 16} (1967), 209-292; English transl. in {\it Trans.
Moscow Math. Soc.,} {\bf 16} (1967).

\bibitem{SkKov}
O.A. Kovaleva and A.L. Skubachevskii,  Solvability of nonlocal
elliptic problems in weighted spaces, {\it Mat. Zametki,} {\bf 67}
(2000), 882-898; English transl. in {\it Math. Notes,} {\bf 67}
(2000).

\bibitem{LM}
J.L. Lions and E. Magenes, {\it Non-homogeneous boundary value
problems and applications, Vol. I}, Springler, Berlin, 1972.
\bibitem{NP}
S.A. Nazarov and B.A. Plamenevskii, {\it Elliptic Problems in
Domains with Piecewise Smooth Boundaries,} Moscow: Nauka, 1991;
English transl. in De Gruyter Expositions in Mathematics, {\bf
13}. Walter de Gruyter Publichers, Berlin -- New York, 1994.

\bibitem{Picone}
M. Picone, Equazione integrale traducente il pi\`u generale
problema lineare per le equazioni differenziali lineari ordinarie
di qualsivoglia ordine, {\it Academia nazionale dei Lincei. Atti
dei convegni,} {\bf 15} (1932), 942-948.

\bibitem{Podiap}
V.V. Podiapolskii, Completeness and basisness by Abel of a system
of root functions of some nonlocal problem, {\it Differentsial'nye
Uravneniya,} {\bf 35} (1999), 568-569; English transl. in {\it
Differential Equations,} {\bf 35} (1999).


\bibitem{SkMs83}
A.L. Skubachevskii, Nonlocal elliptic problems with a parameter,
{\it Mat. Sb.,} {\bf 121} ({\bf 163}) (1983), 201-210; English
transl. in {\it Math. USSR Sb.,} {\bf 49} (1984).

\bibitem{SkMs86}
A.L. Skubachevskii, Elliptic problems with nonlocal conditions
near the boundary, {\it Mat. Sb.,} {\bf 129} ({\bf 171}) (1986),
279-302; English transl. in {\it Math. USSR Sb.,} {\bf 57} (1987).

\bibitem{SkDu90}
A.L. Skubachevskii, Model nonlocal problems for elliptic equations
in dihedral angles, {\it Differentsial'nye Uravneniya,} {\bf 26}
(1990), 120-131; English transl. in {\it Differential Equations,}
{\bf 26} (1990).

\bibitem{SkDu91}
A.L. Skubachevskii, Truncation--function method in the theory of
nonlocal problems, {\it Differentsial'nye Uravneniya,} {\bf 27}
(1991), 128-139; English transl. in {\it Differential Equations,}
{\bf 27} (1991).

\bibitem{SkJMAA}
A.L. Skubachevskii, On the stability of index of nonlocal elliptic
problems, {\it Journal of Mathematical Analysis and Applications,}
{\bf 160} (1991), 323-341.

\bibitem{SkBook}
A.L. Skubachevskii, {\it Elliptic Functional Differential
Equations and Applications,} Basel--Boston--Berlin, Birkh\"auser,
1997.

\bibitem{SkRJMP}
A.L. Skubachevskii, Regularity of solutions for some nonlocal
elliptic problem, {\it Russian J. of Mathematical Physics,} {\bf
8} (2001), 365-374.


\bibitem{Sommerfeld}
A. Sommerfeld, Ein Beitrag zur hydrodinamischen Erkl\"arung der
turbulenten Flussigkeitsbewegungen, {\it Proc. Intern. Congr.
Math., Rome, 1908, Reale Accad. Lincei. Roma.,} {\bf 3} (1909),
116-124.

\bibitem{Tamarkin}
J.D. Tamarkin, {\it Some General Problems of the Theory of
Ordinary Linear Differential Equations and Expansion of an
Arbitrary Function in Series of Fundamental Functions}, Petrograd,
1917; adridged English transl. in {\it Math. Z.,} {\bf 27} (1928),
1-54.

\bibitem{Ventz}
A.D. Ventsel', On boundary conditions for multidimensional
diffusion processes, {\it Teoriya Veroyatn. i ee Primen.,} {\bf 4}
(1959), 172-185; English transl. in {\it Theory Prob. and its
Appl.,} {\bf 4} (1959).

\bibitem{Vishik}
M.I. Vishik, On general boundary value problems for elliptic
differential equations, {\it Trudy Moskov. Mat. Obshch.,} {\bf 1}
(1952), 187-246; English transl. in Amer. Math. Soc. Transl. (2),
{\bf 24} (1963).


\end{thebibliography}
\end{document}